\documentclass[11pt,a4paper]{article}
\usepackage{t1enc}
\usepackage[latin1]{inputenc}
\usepackage[english]{babel}
\usepackage{natbib}
\usepackage{url}
\usepackage{amssymb,stmaryrd}
\usepackage{caption}
\usepackage{graphicx}

\usepackage[
colorlinks=true,backref,pagebackref]{hyperref}

\newcommand{\beq}{\begin{equation}}
\newcommand{\eeq}{\end{equation}}
\newcommand{\beqarr}{\begin{eqnarray}}
\newcommand{\eeqarr}{\end{eqnarray}}
\newcommand{\beqa}{\begin{eqnarray*}}
\newcommand{\eeqa}{\end{eqnarray*}}
\begin{document}
\thispagestyle{empty}

\title{\bf E. Cartan's  attempt at  bridge-building between Einstein and the Cosserats -- or how translational curvature  became to be known as {\em torsion}}
\author{Erhard Scholz\footnote{University of Wuppertal, Faculty of  Math./Natural Sciences, and Interdisciplinary Centre for History and Philosophy of Science, \quad  scholz@math.uni-wuppertal.de}}
\date{09. 10. 2018 } 
\maketitle
\begin{abstract}
\'Elie Cartan's ``g\'en\'eralisation de la notion de courbure'' (1922) arose from a creative evaluation of the geometrical structures underlying both, Einstein's theory of gravity and the Cosserat brothers generalized theory of  elasticity. In both theories groups operating in the infinitesimal  played a crucial role.  To judge from his publications in 1922--24, Cartan developed his concept of generalized spaces  with the  dual context of general relativity and non-standard elasticity in mind. In this context it seemed natural to express the translational curvature of his new spaces by a rotational quantity (via a kind of Grassmann dualization). So Cartan  called his translational curvature "torsion" and  coupled it to a hypothetical rotational momentum of matter several years before spin was encountered in quantum mechanics. 
\end{abstract}

\setcounter{tocdepth}{2}
{\small \tableofcontents}

\section*{Introduction}
In a series of  notes in the {\em Comptes Rendus} of the Paris Academy of Sciences submitted between February and April 1922  \'Elie Cartan  sketched the basic ideas of a new type of geometry which was centrally based on employing the method of differential forms as far as possible \citep{Cartan:1922a,Cartan:1922[57],Cartan:1922[58],Cartan:1922[59],Cartan:1922[60],Cartan:1922[61]}. The notes were an outgrowth of his investigations of  Einstein's gravity theory under the perspective of his own ideas  in  differential geometry. He completed this first round of publications  by  applying  his methods to the  problem of space as it had  recently been re-formulated by Hermann Weyl in the light of relativity.   A detailed presentation of the ideas followed during the next years.\footnote{See \citep{Chorlay:Cartan,Nabonnand:2016Cartan}.}

One aspect of Cartan's peculiar approach to differential geometry consisted in  formulating the curvature concept of   Riemannian differential geometry  in terms of differential forms with values in the  inhomogeneous Euclidean group operating in the infinitesimal neighbourhoods of any point. But the core of his new geometry lay elsewhere; it  generalized the concept of curvature in two respects.
The first generalization consisted in adding a translational component to the connection  and, correspondingly,  to the curvature. For the latter  he chose the somehow surprising name ``torsion''. The earliest public presentation of his idea was given in his note  {\em Sur une g\'en\'eralisation de la notion de courbure de Riemann et les espaces \`a torsion} of February, 22nd,  \citep{Cartan:1922[58]}.

 The second, perhaps even more consequential, generalization lay in his proposal for allowing  different types of groups operating in the infinitesimal neighbourhoods, rather than just concentrating on the group of Euclidean motions (respectively their Lorentzian counterpart, the Poincar\'e group). This made it possible to study  various types of geometries arising from  the  conformal, the affine, the projective groups, or even more general Lie groups, with their respective pairing of inhomogeneous/homogeneous constituents.
 By this move  Cartan  reshaped the Kleinian program of structuring different types of geometry according to their automorphism groups in the context of differential geometry.  In his note of March 13th, {\em Sur les espaces g\'en\'eralis\'ees et la th\'eorie de la Relativit\'e} \citep{Cartan:1922[59]} the idea was first stated in  some generality.  The double aspect of infinitesimalizing the Kleinian view of geometry and of taking into account a translational component of connection and curvature was crucial for Cartan's {\em espaces g\'en\'ealis\'ees} which later became to known as {\em Cartan spaces}.\footnote{For a modern presentation see \citep{Sharpe:DiffGeo}.}

The following paper  concentrates on the first aspect of Cartan's generalization of differential geometry  and the peculiar contexts which lay at the base of the, prima facie paradoxical,  terminology of {\em torsion} for the {\em translational} component of the curvature. In the paper in which Cartan announced this new concept he described it in quite intuitive terms. He expressed the difference of his approach to classical (Euclidean) geometry similarly to what had been done by Levi-Civita and Weyl. That is, he considered the change a vector would undergo, if it is transported along an infinitesimal closed path according to the rules established by the  generalized connection:  
\begin{quote}
En d\'efinitive, \`a tout contour ferm\'e infiniment petit de l'espace donn\'e sont associ\'ees une translation et une rotation infiniment petites (. . . ) qui manifestent la divergence entre cet espace et l'espace Euclidien \citep[p. 594]{Cartan:1922[58]}.\footnote{{\em Definitely, to any infinitesimally closed curve of the space an infinitesimal translation and a rotation are associated (\ldots); they express the divergence between this space and Euclidean space.} -- Translations in emphasized letters by ES; other translations in quotes with source indicated.}
\end{quote}
The mentioned infinitesimal translation and rotation expresses the curvature properties of the space. 
Cartan immediately identified the well known case of Riemannian geometry with its Levi-Civita connection as the situation in which the translational component of the curvature vanishes. 

A little later in the note he came back to the difference to the more classical geometries again and introduced a new terminology for the translational curvature mentioned  above:
\begin{quote}
Dans les cas g\'en\'eral o\`u il y a une translation associ\'ee \`a tout contour ferm\'e infiniment petit, on peut dire que l'espace donn\'e se diff\'erencie de l'espace euclidien de deux mani\`eres: 1$^{\circ}$ par une {\em courbure}  au sens de Riemann, qui se traduit par la rotation; $2^{\circ}$ par une {\em torsion}, qui se traduit par la translation \citep[594f.]{Cartan:1922[58]}.\footnote{ {\em To any closed infinitesimal loop there is generally an associated translation; in this case  one can say that the given space  differs from Euclidean space in two respects: 1. by a \underline{curvature} in the sense of Riemann, which is expressed by the rotation; 2. by a \underline{torsion} which is expressed by the translation} (emphasis in the original). }
\end{quote}
But why did he  call the translational curvature  ``torsion''?

A first clue follows immediately; but at first glance it enhances the riddle and introduces an even wider {\em quid pro quo}: \label{quid pro quo}
\begin{quote}
La rotation peut \^etre repr\'esent\'ee par un vecteur d'origine $A$ et la translation par un couple (ibid.).\footnote{\em The rotation can be represented by a vector of origin $A$ and the translation by a couple.}
\end{quote}
Now everything has been turned upside down, rotations were expressed by  vectors, translations by  couples. 

The last word of the sentence indicates that a mechanical context stood behind this move. In fact, Cartan 
indicated that one can study the equilibrium of an elastic medium in terms of his connection and curvature. This  led him to formulate  a geometrical picture of the constellation of forces:
\begin{quote}
On a ainsi une image g\'eom\'etrique d'un
milieu mat\'eriel continu en \'equilibre, mais dans le cas o\`u ces forces
se manifesteraient sur chaque \'el\'ement de surface, non seulement
par une force unique (tension ou pression), mais par un couple
(torsion)  \citep[594]{Cartan:1922[58]}.\footnote{\em One thus has  a geometrical picture of a  continuous material medium in equilibrium, but in the case where the forces 
express themselves not only by a single force (tension or pressure) but also by a couple (torsion).}
\end{quote}

By {\em couple} Cartan referred to  the traditional (18th and 19th century) expression for a rotational momentum (torque) by a pair of forces of the same norm, acting along different parallel lines in opposing orientations. This might superficially explain the rephrasing of translational curvature as ``torsion''. But  for the unprepared reader it still  remains a riddle why Cartan identified the  infinitesimal {\em rotations} with forces (vectors) and infinitesimal {\em translations} with rotational momenta (couples). From a purely geometrical point of view this identification would not  appear particularly plausible.
But at the end of the note Cartan gave a hint for the motivation of such an interchange. He indicated that
\begin{quote}
  \ldots les  consid\'erations pr\'ec\'edentes (\ldots) du point de vue m\'ecanique, 
s'apparentent aux beaux travaux de MM. E. et F. Cosserat sur, l'action Euclidienne \ldots''. (ibid.)\footnote{\em  \ldots from the mecanical point of view (\ldots) the preceding considerations are similar to the beautiful works of the Messieurs E. and F. Cosserat on the Euclidean action \ldots }
\end{quote} 
In addition he mentioned  another link, namely to  to H. Weyl's  studies of the problem of space; but this does not lead us further for our question.\footnote{For Weyl's space problem see, among others, \citep{Bernard:2018Paris,Scholz:2016Weyl/Cartan}; for the elasticity of the Cosserats \citep{Brocato/Chatzis:Cosserat,Pommaret:1997,Hehl/Obukhov:2007}.}

If we want to understand  the background of Cartan's choice of  terminology for the translational curvature we  have to reconstruct the historical context of the { unconventional theory of elastic media} of the brothers \'Eug\`ene and Fran\c{c}ois Cosserat, which Cartan referred to. On the other hand, the geometrical picture of the elastic medium Cartan had in mind arose from his way of reading Einstein's gravity theory in a mathematical analogy to elasticity. In order to understand Cartan's intentions epressed in the note \citep{Cartan:1922[58]}  we have to follow the traces of a ``threefold knot''  tied by  Cartan between the {\em mathematical methods} developed for his new type of geometry,  {\em Einstein's theory} of gravity, and the {\em generalized theory of elasticity} of the Cosserats.

We therefore start this paper with a short description of Cartan's mathematical arsenal used for constituting his generalized geometries (section 1), continue with a r\'esum\'e of  Cosserat elasticity and its historical context (section 2) before we shed a glance at Cartan's reading of Einstein gravity (section 3). This allows us to reconstruct how Cartan linked these three components  in an intriguing interplay between  his geometrical picture of a Cosserat type elasticity theory and a (speculative) generalization of Einstein gravity by torsion (section 4). We then look back at his  practice of organizing the three-sided  interplay between mathematics/geometry, elasticity theory, and gravity (section 5), and give some indications of its repercussions  on the work of physicists in the second half of the 20th century (section 6).

\section{\small A  short outline of  basic ideas of Cartan geometry \label{section Cartan geometry}}

 The usual differential geometric description of a {\em metric}  $ds^2$ (Euclidean, Minkowski or (pseudo-)Riemannian) uses the differentials $dx_i$ of the coordinates of a point  $x=(x_1, \ldots, x_n)$
	\[  ds^2 = \sum_{i=1}^n g_{ij} \ dx_i^2 dx_j^2 \] 
This corresponds to a choice of a coordinate basis in the tangent spaces (infinitesimal neighbourhood) of any point \ldots \\
	\hspace*{4cm}{\includegraphics[angle=0,scale=0.2]{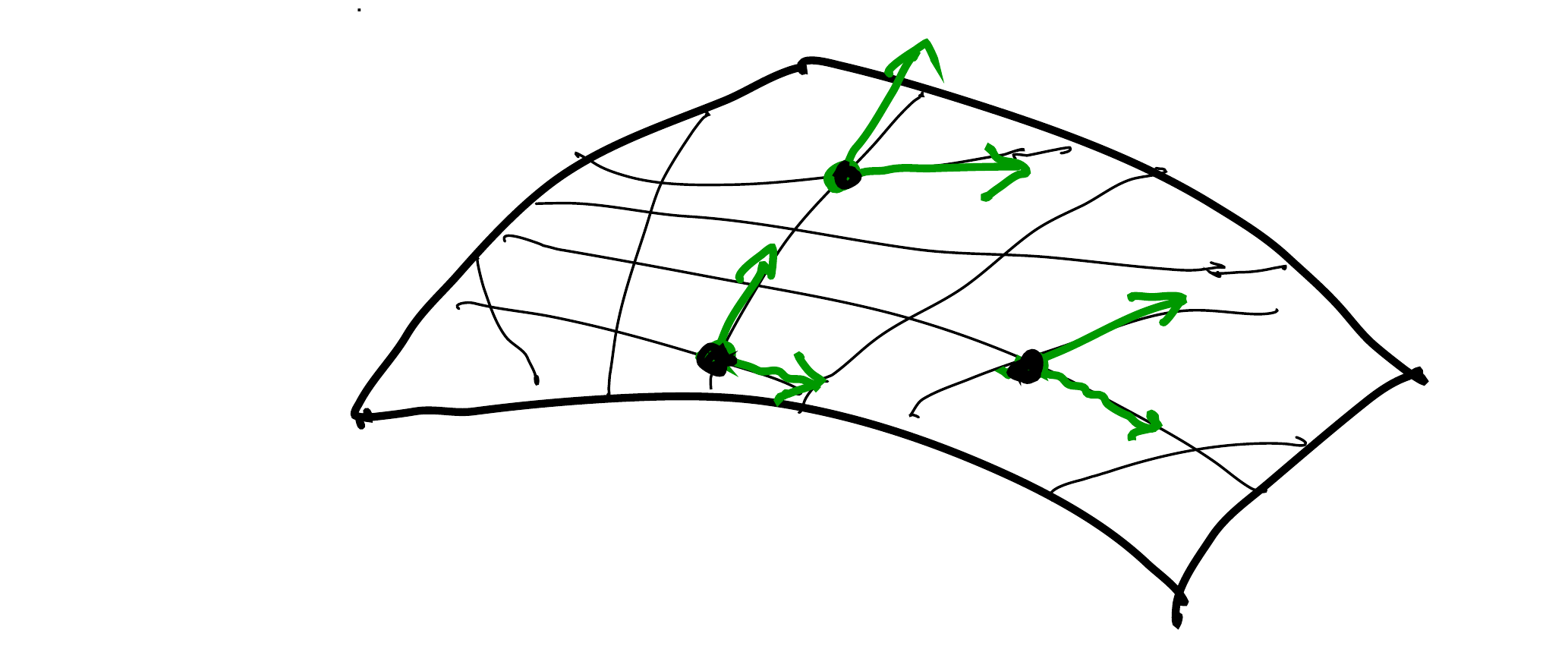} } \\
and an expression of the metric with regard to this basis.

Cartan, in contrast, preferred to	 describe the metric in terms of  {differential forms  $\omega_1, \ldots \omega_n$}  which diagonalize the metric:
\\[-0.5em]	
\beq ds^2 = \sum_{i=1}^n \epsilon_i \, \omega_i^2, \qquad   \epsilon_i = \pm 1
\label{metric}	\eeq
 The  $\epsilon_i$ (used by Cartan himself) account for different signatures of the metric, most importantly Euclidean/Riemannian and Minkowski/Lorentzian.\footnote{\citep[p. 150, eq. (9)]{Cartan:1922a}} 
 
 This form can be arrived at by  linear algebraic considerations in each infinitesimal neighbourhood.
In his papers of 1922ff. Cartan emphasized that geometrically the diagonalization indicates  a choice of point-dependent (``mobile'') {\em orthonormal  reference systems}.\footnote{The paper \citep{Cartan:1922a} was written in 1921 and published only in the following year. Cartan remarked that the ``germs'' of his new geometry can be found at the beginning and the end of this paper. Before it was published, Cartan announced the basic ideas of his new geometry in several {\em Comptes Rendus} notes, \citep{Cartan:1922[57],Cartan:1922[58],Cartan:1922[59],Cartan:1922[60],Cartan:1922[61],Cartan:1922[62]}. Technical details followed in his long m\'emoire {\em Les vari\'et\'es a connexion affine et la th\'eorie de la relativit\'e g\'en\'eralis\'ee} \citep{Cartan:1923/24,Cartan:1925} and subsequent publications.}
  Cartan called them Euclidean reference systems,  ``syst\`eme de r\'eference euclidien'' (ibid., p. 151) or  ``tri\`edre trirectangulaire'' \citep{Cartan:1922[58]}  etc. and denoted them 
 by  { $e_1, e_2, \ldots, e_n$} (orthonormal basis ONB, or frame). \\
\hspace*{2cm} \includegraphics[angle=0,scale=0.4]{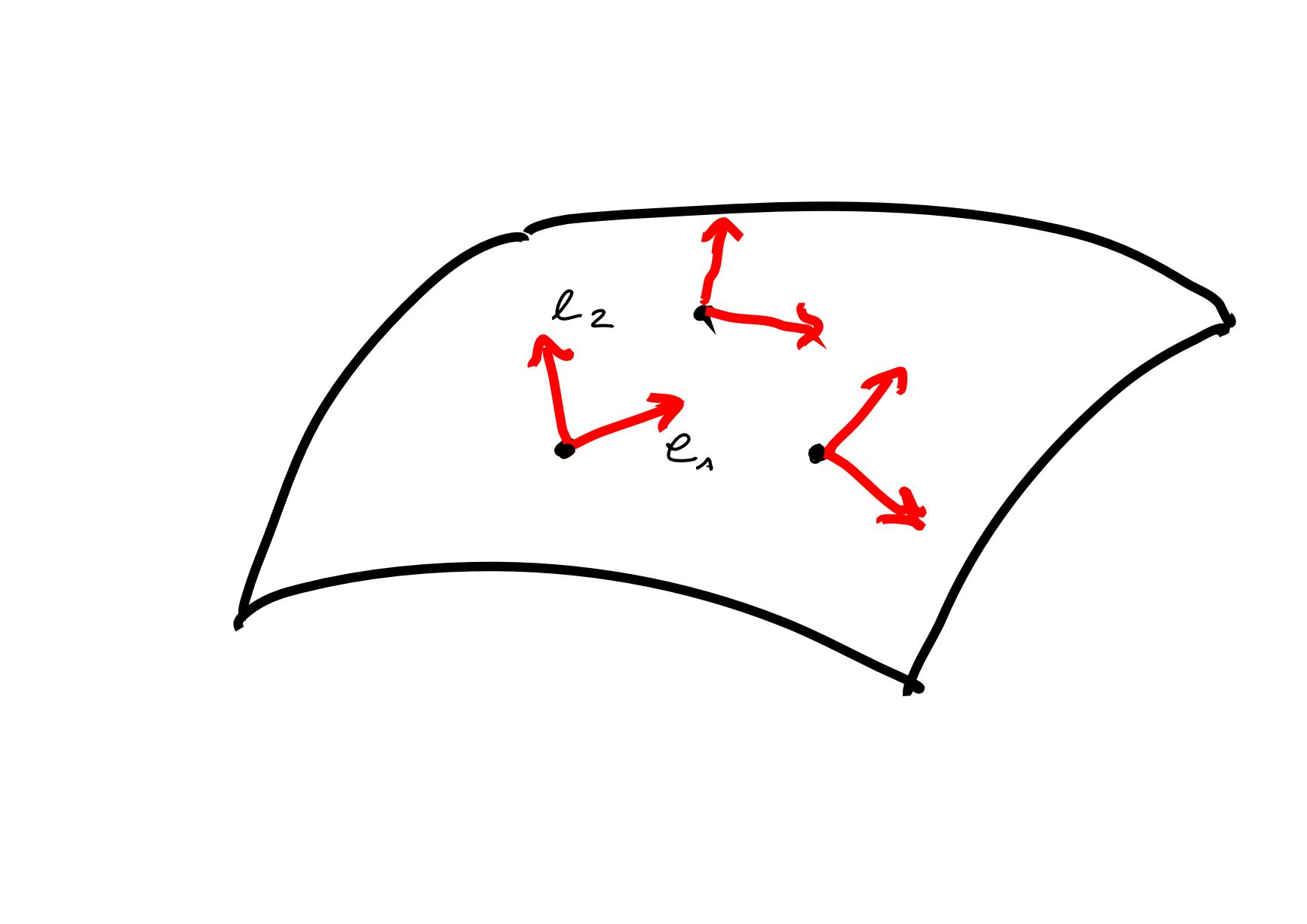} \\
If we introduce the analogous symbols for the  dual basis (at evey point),
\beq \{e^1, \ldots e^n\}  \quad \mbox{dual basis of 1-forms to ONB} \quad \{e_1, \ldots, e_n\}, , \label{ONB-dual}
\eeq 
Cartan's $\omega_i$ turn out to be nothing but these, $\omega_i = e^i$. In a coherent use of lower and upper indexes one therefore better writes Cartan's component forms as $\omega^i$. In fact  Cartan often, although not always, used upper and lower indices like in the tensor calculus,\footnote{Upper indices ones for vector like, and lower ones for differential form like transformation behaviour under change of coordinates or reference systems.}
 e.g.  $\omega^i_{\; k}$ in place of $\omega_{ik}$.
He also applied the Einstein summation convention abbreviating, e.g.,   $ \sum_k \omega_{ik}\omega_{kj}$ by $ \omega^i_{\; k}\omega^k_{\; j}$ etc.

If one moves  between infinitesimally close points $x, x'$ the reference systems undergo an   infinitesimal rotation given by a system of coefficients $(\omega_{ij})$ depending on the start point $x$ and $\delta x= x'-x$:  \\[0.8em]
\hspace*{2cm} \includegraphics[angle=0,scale=0.5]{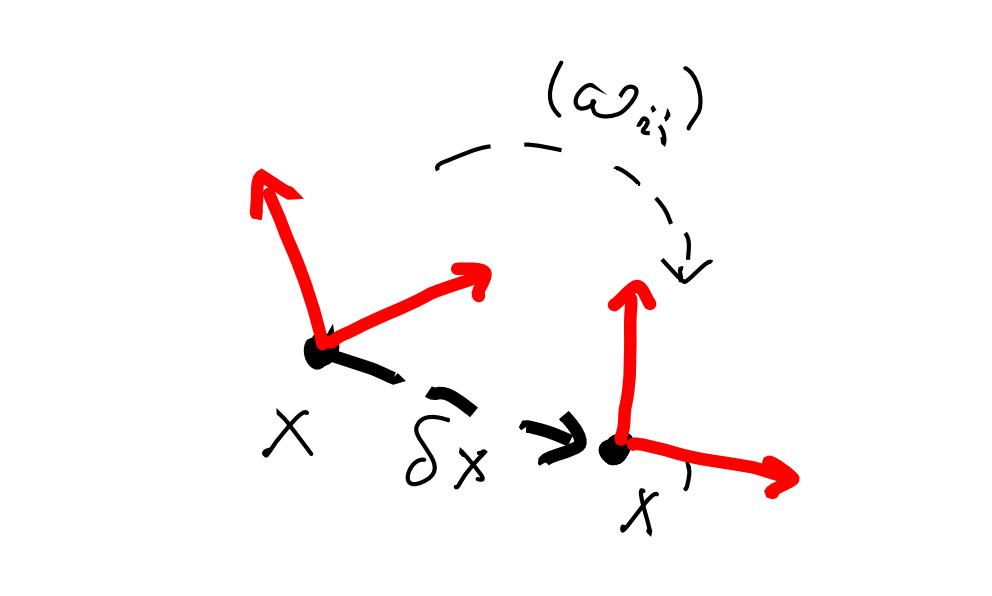} \\
Cartan realized that the coefficients of {$\omega_{i j}$} can be understood as a system of differential forms   (antisymmetric in the indices $i,j$). They encode the  {\em rotational connection} of the space.

By analogy Cartan interpreted the    {$\omega_i$} as  assigning to any $\delta x$ a  translational shift of the reference system identical to $\delta x$:\\  
\hspace*{2cm} \includegraphics[angle=0,scale=0.12]{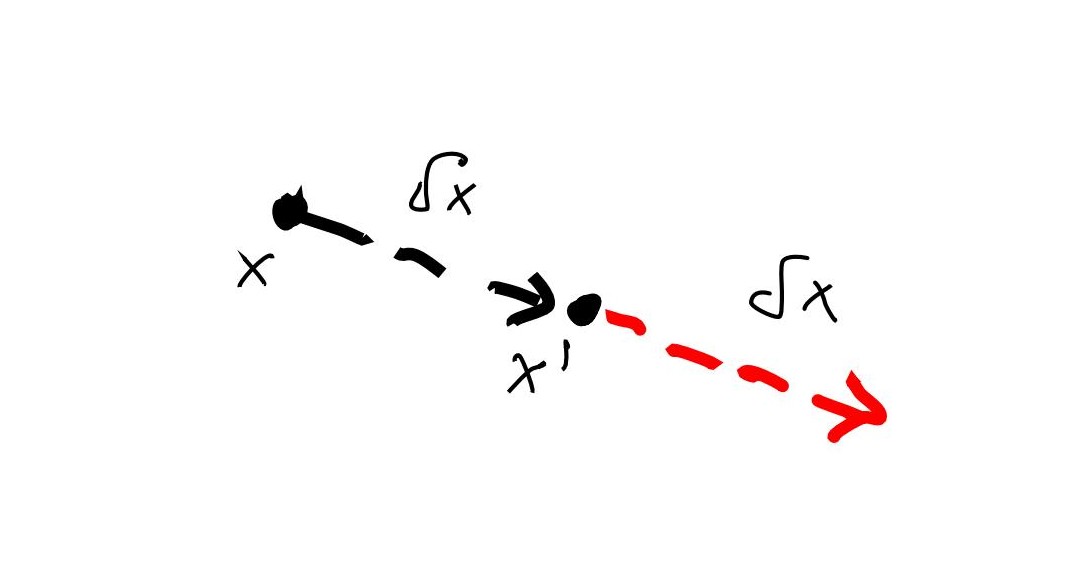}  \\
This was  a new idea which paved the way for the first of Cartan's two innovations mentioned above. In addition to the role of the $\omega_i$ for representing the metric in diagonal form (and for specifying a ``tri\`edre trirectangulaire'') he used them for  assigning a translation with components $\omega_i(\delta x)$ to any infinitesimal  shift $\delta x$ from the point $x$ to an infinitesimally close one $x'$ \citep[p. 152]{Cartan:1922a}. This was a first step towards turning the $\omega_i$ into a {\em translational connection} which complements  the rotational one of the reference systems:\\ 
\hspace*{2cm} \includegraphics[angle=0,scale=0.12]{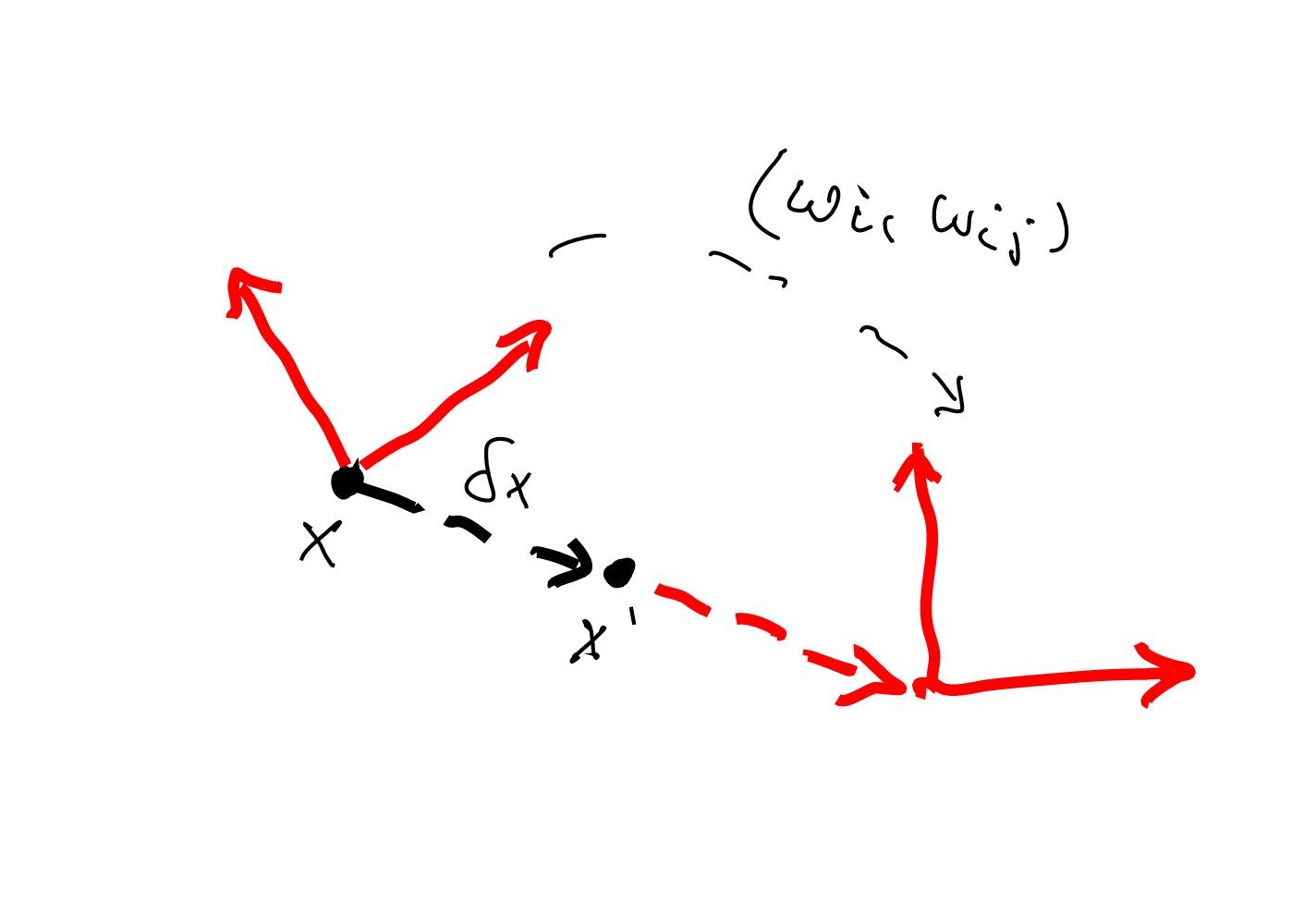} \\

An important feature of Cartan's approach was that both parts of connection, the rotational and the translational one, were given  component-wise by (real valued) differential forms. Present day readers may prefer to read them more collectively as two  differential forms, one with values in the Lie algebra of the rotational group $\overline{\omega}=(\omega^i_{\;j})$ and one with values in the translations $\omega=(\omega^i)$. Moreover, a present reader might like to see an explicit expression for the covariant derivative $\nabla$ of Cartan's connection $\widetilde{\omega}=(\omega, \overline{\omega})$, which would generalize the Levi-Civita connection of the metric (\ref{metric}).\footnote{For a lucid modern presentation see \citep{Sternberg:Curvature}; for more technicalities \citep[appendix]{Gasperini:2017}.} 
Cartan   emphasized calculations which could be expressed in the calculus of differential forms, rather than rewriting the bulk of Ricci's and Levi-Civita's covariant tensor calculus in his symbolism. 

For calculating the analogue of exterior differentials of the connection forms 
Cartan had to take the rotational coefficients into account. If both generalized exterior differentials  are zero, 
\beqarr d \omega_i + \epsilon_i \sum_k \omega_{ik} \omega_k &=& 0,  \label{flatness condition} \\
\nonumber d \omega_{ij} + \sum_k \epsilon_k \omega_{ik}\omega_{kj} &=& 0,
\eeqarr
so Cartan noticed, the space is {\em  Euclidean} (flat).\footnote{\citep[pp, 145, 148]{Cartan:1922a}}

 But in general this need not be the case,  and  one encounters an  {\em espace g\'en\'eralis\'e} (generalized space).
  If one then lets a point $M$ traverse an infinitesimal loop starting and ending in $A$ 
  \begin{quote}
\ldots  {\em  on ne retrouvera pas dans l'espace euclidien le tri\`edre initial} mais il faudra, pour l'obtenir, effectuer un d\'eplacemant compl\'ementaire dont les composants sont bien d\'efinies par rapport au tri\`edre initial \citep[p. 594, emphasis in the original]{Cartan:1922[58]}.\footnote{\em \ldots \underline{one does not find the initial three-frame in the Euclidean space} ({\em meant is the tangent space in modern terms, ES}),  in order to arrive at it one rather has to apply a complementary displacement the components of which are well defined with regard to the original three-frame.} 
  \end{quote}
Cartan noted explicitly that this   {\em d\'eplacement complementaire} is {\em independent} of the choice of reference systems. Today it is  called the {\em Cartan curvature} of the space.

The translational and rotational {\em d\'eplacement complementaire} $\Omega^i$ and  $ \Omega^i_{\;j}$, i.e. the deviations from zero of the above given expressions are  the  2-forms:
	\beqarr  \Omega^i &=& d \omega^i + \sum_k \omega^i_{\; k} \omega^k =A^i_{\;j k} \, [\omega^i \omega^j ] \qquad \mbox{( {``torsion''})} \label{equ structure}\\
	 \Omega^i_{\;j}  &=&d \omega^i_{\;j} + \sum_k \omega^i_{\; k}\omega^k_j = A^i_{\;jkl}\, [\omega^k \omega^l ]  \qquad \mbox{( {``courbature''}),}  \label{equation courbature}
	\eeqarr
where the square brackets denote alternating products. They  were adapted  by  Cartan from Grassmann.\footnote{Later authors, more precisely E. K\"ahler, introduced the now common symbolism $[a b] = a \wedge b$.} 
Cartan called the equations (\ref{flatness condition}) and (\ref{equ structure}) ``les \'equations de structure'' (structural equationis) of the generalized space \citep[p. 368]{Cartan:1923/24}. For vanishing torsion the $A^i_{\;jkl}$ characterize the Riemannian curvature in Cartan's symbolism \citep[p. 154]{Cartan:1922a}. But it may also happen that the rotational curvature vanishes, while the torsion is non-trivial.

A curve with a  tangent vector field which is parallel in the sense of the Cartan connection is called an  {\em autoparallel}, while a curve of extremal length is called a {\em geodesic}. In Riemannian geometry both concepts agree, but in Cartan geometry (modeled on the Euclidean or pseudo-Euclidean group) they usually fall apart.  
But this need not be necessarily so.  Cartan gave a simple example of  a structure in dimension $n=3$ with vanishing rotational curvature and non-trivial torsion, in which geodesics and autoparallels coincide \citep[p. 595]{Cartan:1922[58]}.\footnote{In the context of the studies of (generalized) Cosserat media  this structure   attracts attention until today  as an example with intriguing geometrical properties \citep{Lazar/Hehl,Hehl/Obukhov:2007}.}
 
Cartan generalized this approach to allow for more general groups than the orthogonal ones, operating in the infinitesimal neighbourhoods. At the moment we need not follow this generalization in more details;  but in general, the metric lost its central place and and the ``tri\`edre trirectangulaires'' had to be replaced by more general ``r\' eperes''. Cartan called the arising spaces {\em espaces non-holonomes} (non-holonomous spaces).\footnote{For the historical background of this terminology see \citep{Nabonnand:Cartan_2009}.} 
During the 1920s  he  studied such spaces of  increasingly complex type with the following groups: 
 \begin{itemize}
\item The Poincar\'e group in papers on the geometrical foundation of  general relativity \citep{Cartan:1922[59],Cartan:1923PoS}, \citep{Cartan:1923/24}. For torsion  $\Omega ^i =0$ such a Cartan space reduces to a Lorentzian manifold. Cartan could use this reduction for treating Einstein's theory in his own geometric  terms.
\item The inhomogeneous similarity group. For torsion $ =0$, this case reduces to Weylian manifolds    \citep{Cartan:1923PoS}.
\item The conformal group \cite{Cartan:1922[60]}.
\item The projective group \cite{Cartan:1924[70]}.
\end{itemize}

 In this way, Cartan developed a wide  conceptual frame for studying different types of differential geometries, Riemannian, Lorentzian, Weylian, affine, conformal, projective. All were  enriched  by the possibility to allow for the new phenomenon of torsion, and all  arose from Cartan's unified method of adapting the Kleinian viewpoint to infinitesimal geometry. But if we want to understand his first papers of the year 1922 and the immediately following ones, we have to know a bit the ``beaux travaux de M. E. et F. Cosserat''  \citep{Cartan:1922[58]}.

\section{\small Generalized elasticity theory  \label{section elasticity}}
In the early 19th century a group of mainly French authors developed the foundations  of the linear elasticity theory of solid bodies. A.J. Fresnel (1821) and C.L. Navier (1827) derived their theories on the basis of a molecular theory of matter with central forces acting between the discrete units of matter. When A.-L. Cauchy jumped in between 1823 and 1828, he first approached the question from the point of view of a continuum theory of matter and derived his influential representation of the linear relationship between the {\em strain} matrix  characterizing the deformation of the material and the {\em stress} matrix (both later understood as tensors) from a phenomenological {\em Ansatz}.\footnote{For Cauchy's contributions to elasticity see \citep{Dahan:Cauchy,Belhoste:Cauchy}.}

 But the molecular theory of matter behind these different approaches remained dominant. In 1827 also Cauchy presented a derivation of Navier's equations  on the basis of a molecular approach. A year later S.D. Poisson developed the linear elasticity theory of molecular matter a step further and brought it home to the Laplacian program of physics, which in the meantime had come under attack from different sides (Fourier's theory of heat, electricity, magnetism, optics) \citep{Fox:LaplacianPhysics}.

Poisson's theory was built upon the hypothesis of central forces acting between point-like centers inside the radius of a ``molecular sphere'' outside of which the forces are no longer to be felt. The phenomenological forces in the material on a surface element were derived by summing up all the forces in the range of the ``molecular spheres'' of points intersecting or touching the  surface element. For isotropic solid matter the calculations resulted in a linear relation between strain (deformation) and stress (surface forces), which depended on a single material constant. The basic structure of the theory seemed empirically convincing; but  with increasing precision of experimental techniques between the 1840s and 1870s  the 1-parameter assumption turned out to be untenable even for isotropic matter;  a second elastic constant had to be assumed to fit the data. Even worse, around the middle of the century the assumption of pointlike molecules of the Laplacian program became undermined also from another side: The improvements in theories of crystal structure, in particular  A. Bravais' theory of crystal matter (late 1840s to early 1860s)  indicated that directional aspects might well play a role also for the elastic properties of  matter.\footnote{For more details on this development see
\citep{Capecchi_ea:2010,Fox:LaplacianPhysics,Timoshenko:1953}.}  

An alternative approach to the theory of elasticity was proposed by George Green in 1838. He  avoided any hypothesis about the basically unknown molecular structure of matter and based his analysis on a potential function $\phi$ from which the forces in the elastic medium could be derived by very general formal considerations. Although this approach  led to quite acceptable results, including the empirically necessary two elastic constants in the case of an isotropic medium, Green's theory did not manage to replace the research program following the  molecular hypothesis  \citep[pp. 217ff.]{Timoshenko:1953}. But it became an important input for the generalized theory of elasticity of the brothers Cosserat to whom Cartan referred in his note of 1922. 

In the late 1880s  {\em  Woldemar Voigt} (1850--1919) gave a detailed analysis of the actual status of the molecular theory of elasticity in a report to the {\em G\"ottingen Gesellschaft der Wissenschaften}  \citep{Voigt:1887}. He carefully reviewed the molecular elasticity theory of the French tradition and 
 proposed a refinement of it, which would take into account that the molecules are extended bodies of different shapes. In general, the form of the molecules  breaks the rotational symmetry of the old pointlike force centers; thus not only the coordinates of the centers of the molecules, but also their directional properties had to be considered. 
 
 As a result, the molecular interactions could no longer be represented by forces alone but had to be complemented by the consideration of  rotational momenta, torque, which depend on the relative ``polarity'' of the molecules.\footnote{``Wir denken uns das homogene krystallinische Medium bestehend aus einem System von Molek\"ulen, welche durch ihre Wechselwirkungen einander im Gleichgewicht halten. Diese Wechselwirkungen sind Kr\"afte und Drehmomente, deren Componenten in unbekannter Weise mit der relativen Lage der Molek\"ule variieren.'' \citep[p. 5]{Voigt:1887}\\
 {\em We conceive the homogeneous crystalline medium as consisting of a system of molecules which stand in equilibrium by their mutual interactions. These interactions are forces and rotational momenta, the components of which vary  with the relative position of the molecules in an unknown way. }}
For a full representation of the position and the orientation of the  molecules the coordinates of their barycenters and the directions of a system of axis, tied to the molecule and changing from one to the other had  to be taken into account.\footnote{``Da die Molek\"ule nach unserer Annahme eine Polarit\"at besitzen, so muss man sie wie  endliche K\"orper behandeln und ihre Lage ausser durch die Coordinaten ihres Schwerpunktes noch durch die Richtung eines fest mit ihnen verbundenen Axensystemes bestimmen.'' (ibid., p. 6)\\
{\em According to our assumption   the molecules  possess a polarity,  one therefore has to treat them like finite bodies and has to specify their position  in addition to the coordinates of their barycenter by the direction of an axis system rigidly tied to them.}
}
From a mathematical point of view, Voigt's description resembled point dependent {\em r\'ep\`eres mobiles} linked to the  different orientations (``polarisations'') of the molecules in a material structure. But neither he nor mathematicians at the time took up this analogy.

 For studying   equilibrium conditions on the macro-level, Voigt  considered forces and rotational momenta  on surface or volume elements, given with regard to an axis system by the components $(Y,Y,Z)$ and $(L,M,N)$ respectively. They  came about from the summation of the corresponding actions on the micro-level and  had to be studied in  the rest state and, if subject to external forces, in a deformed state (ibid. p. 10). A clear and quite detailed study of Voigt's further derivation is given in \citep{Capecchi_ea:2010}. We need not go into  the  details here, because in the course of his calculations Voigt introduced the assumption that for all practical purposes {\em the point dependence of the rotational  momenta} induced in the material even by deformations {\em could be neglected}.\footnote{ ``\ldots sind die in den Ausdr\"ucken f\"ur die Drehungsmomente vorkommenden Coeffizienten als unendlich klein gegen die in den Componenten $X_x\ldots$ auftretenden anzusehen. Dies hat den Effekt ihre Differentialquotienten neben den \"ubrigen Gliedern zu vernachl\"assigen sind, --- in \"Ubereinstimmung mit der Umstande, dass bei allen bekannten Problemen an der Oberfl\"ache der elastischen K\"orper $\overline{L}_n, \overline{M}_n, \overline{N}_n$  (the rotational momentum represented as a vector normal to the surface, ES) gleich Null zu setzen ist  \ldots'' \citep[p. 23]{Voigt:1887}. \\
 {\em \ldots the expressions for the coefficients in the rotational momenta have to  be considered as infinitely small with respect to those appearing in  the components $X_x$\ldots This leads to the effect that their differential quotients can be neglected in comparison with the other terms -- in accordance with the fact that in all known problems at the surface of elastic bodies  $\overline{L}_n, \overline{M}_n, \overline{N}_n$  (the rotational momentum represented as a vector normal to the surface, ES) may be equated to zero \ldots}
 }
 
So the bulk of Voigt's enriched structure theory on the micro-level (the point dependence of the axis systems linked to the molecules and their deformations) remained without visible consequences, once one turned to  the phenomenological level. One effect remained however.  Voigt's calculations led to introducing a second parameter resulting from a global rotational momentum which was not present in the older molecular theory.  It 
filled the gap which had arisen between the older molecular theory of elasticity and the experimental findings. In the end, this was the main achievement of Voigt's approach. It soon became accepted and shaped the paradigm of linear elasticity theory at the turn to the 20th century.

Other authors explored  alternatives in the framework of continumm mechanics. Particularly important in  our context  was the joint research of {\em Fran\c{c}ois  Cosserat}
(1853 -- 1914)  and his younger brother   {\em \'Eug\`ene Cosserat} (1866 -- 1931) during the two  decades between 1896 and 1914. Fran\c{c}ois was a civil engineer working for the French railroad system. He  studied at the \'Ecole Polytechnique and graduated at the \'Ecole des Ponts et Chauss\'ees and was a highly theoretical mind. \'Eug\`ene  studied mathematics at the \'Ecole Normal Superieur under the guidance of P. Appell, G. Darboux, G. Koenigs  and E. Picard. After his graduation in 1886 and  a few months of teaching at a Lyc\'ee in Rennes he became an assistant astronomer at the Observatory in Toulouse. 
 Parallel to observational work on binary stars he wrote a dissertation in mathematics with a topic in differential geometry. In 1889 he finished his PhD in Paris  under the supervision of Appell, Darboux and Koenigs. 
   In 1896 he succeeded T. Stieltjes as a professor in mathematics at  Toulouse University. Roughly a decade later (1908) he became the director of Toulouse Observatory and professor of astronomy. He was elected corresponding member of the Paris Acad\'emie des Sciences in 1911 and became a full member in 1919 \citep{Levi:Cosserat}. After the early death of his older brother he discontinued work on elasticity theory. The last joint publication of the two  appeared after  Fran\c{c}ois' death.  It was an extended French version of  A. Voss essay on rational mechanics for the {\em 
  Encyclop\'edie des sciences math\'ematiques pure et appliqu\'ees}. It provided  the occasion for  explaining the wider perspective of their research in rational mechanics \citep[sec. 4]{Brocato/Chatzis:Cosserat}. 
   
The two brothers  studied elasticity theory   in a strictly deductive Lagrangian approach to continuum mechanics, while acknowledging that its aim was a rational understanding of inductively generalized empirical knowledge. Their first paper (``premier m\'emoire'')  appeared in the year \'Eug\`ene became a professor of mathematics in Toulouse \citep{Cosserat:1896}. A series of papers followed; but their main result did not become mature before 1909. In this year they  managed to  derive,  on the basis of two principles, a set of generalized equations for the equilibrium of an elastic medium carrying forces and torques.  They presented their new theory  in several variants to the scientific public  \citep{Cosserat:1909Theorie,Cosserat:1909Note,Cosserat:1909Chwolson}.\footnote{For  list of all common papers of \'Eug\`ene and Fran\c{c}ois Cosserat and a discussion see \citep{Brocato/Chatzis:Cosserat}.}

Their  first principle was the {\em invariance} of the action under transformations of the {\em  inhomogeneous Euclidean group} for  elastic continua of dimension $n=1, 2, 3$ (elastic rod, plate, body). In their terminology they worked with  an  ``action euclidienne''. As a second principle they  characterized the   elements of the elastic continuum  by  {\em  point dependent ``tr\i\`edres''} (orthonormal frames) rather than considering elastic deformations of a simple point continuum.
In this form they took up ideas in   elasticity theory  which had started to consider ``polarized'' (directionally oriented) molecules. They  incorporated them into the continuum mechanics framework  and gave them a form which  nicely corresponded to the methods of Darboux style differential geometry.\footnote{In his PdD dissertation \'Eug\`ene had already investigated  infinitesimal circles as space elements,  combining ideas of Pl\"ucker's generalized space elements with Darboux' differential geometry. }

In their last paper they  described a limit idea underlying this approach. In the older approach to elasticity  ordinary geometric space was considered as an adequate mathematical  representation of the  physical medium (``milieu'').  But, according to the Cosserats,  the studies of elasticity, cristallography, electricity and of light made it necessary to consider  a more complex notion of the continuum (``une notion plus complexe du milieu continu'').
\begin{quote}
\ldots This notion is derived in all generality by a passage to a limit from a discontinuous collection of {\em point systems with an arbitrary number of degrees of freedom}.\citep[p. 72, emphasis  ES]{Cosserat:1915} \footnote{``Cette notion se
d\'eduit dans toute sa g\'en\'eralit\'e, par un passage \`a la limite, de celle d'un ensemble discontinu de syst\`emes de points \`a un nombre quelconque de degr\'es de libert\'e'' citation from \citep[p. xx]{Brocato/Chatzis:Cosserat}. With regard to crystallography one may easily recognize the discrete ``point systems'' as an abstract representation of the lattice structure of polyhedral molecules studied by Bravais or the more refined structures of the 1890s according to the research tradition of Fedorow or Schoenflies.}
\end{quote}
This may be read as a late reflection on the motivations which had brought  them to  study the influence of the directionally oriented elements on the equilibrium conditions in all generality (not only under the restrictive  assumptions used by Voigt).\footnote{The Cosserats saw and commented  the relationship between Voigt's and their work. They were clearly aware of their own achievements; see fn   \ref{fn Voigt-Cosserat}.}

 The Cosserats characterized an  element of the undeformed continuum (``\'etat primitif'' or  ``\'etat naturel'')  by the coordinates $x=(x_1,x_2,x_3)$ of a point $p$ with regard to a fixed Euclidean frame $\mathfrak{O}$  and  orthonormal frame (``tri\`edre trirectangle'') $\{e_1,e_2,e_3\}$ attached to the point and specified by a point dependent rotation $o(x)$ with regard to the reference system $\mathfrak{O}$. For the sake of brevity we denote such an oriented continuum element here by $(x, o(x))$.\footnote{The original notation of the Cosserats for the tri\`edre $M_o x_o, M_o y_o,M_o z$ was given by angle cosinus to the fixed reference system \citep[p. 559]{Cosserat:1909Note}.}
 The coordinates $x$ could be changed by a smooth coordinate transformation. A deformed state $(x', o'(x'))$ of the medium, on the other hand, was  described by transformations $x'=f(x)$ and $o'(x)= g(o(x))$ with smooth functions $f, g$. The range of possible infinitesimal deformations was then characterized by the 9 partial derivatives of  the three components of $f$, which we denote collectively by $\partial f$  and 9 partial derivatives of the angle transformations of type $\partial g$.\footnote{Expressed in coordinates of  the fixed reference system $\mathfrak{O}$, the  Cosserats gave  $\partial f$ as $(\xi_i,\eta_i, \zeta_i)$  and $\partial g$  by $(p_i, g_i, r_i)$, where the index $i=1,2,3$ indicates the partial derivative with regard to $x_i$ \citep[pp. 559, 596]{Cosserat:1909Note}, similar in \citep{Cosserat:1909Theorie}.}
 
They assumed a time-independent action  density $W$ for the deformation of the continuous medium and analyzed it step by step for the dimensions $n=1,2,3$. For dimension $n=3$ the action was of the general form
\beq  \mathfrak{S}= \int_{A_o} W(x,\partial f, \partial g  )\,dx \; , \label{action}
\eeq
 with $A_o$ the space region occupied by the elastic body in the natural state.  Thus their action depended on 21=3+9+9 continuous parameters \citep[p. 559]{Cosserat:1909Note}. Its form was  constrained only by the {\em postulate of invariance under (infinitesimal) Euclidean motions}. Analyzing the variation of  a not necessarily ``natural'', i.e. force free, state they were able to  derive formal expressions for  the  {\em external forces} $(X_i)$ and {\em rotational moments} $(L_i)$ acting on the  volume elements.\footnote{Volume forces may result from  fields  permeating the continuum; but the authors did not discuss the origin of them.} 
 The 
 \beqa
 \mbox{surface densities of} & &  \mbox{force} \;  (F_1,F_2,F_3) \; \quad  \mbox{and torque} \; (J_1,J_2,J_3)\; \\
\mbox{and volume densities of} & &  \mbox{force} \;  (X_1,X_2,X_3)  \quad \mbox{and torque} \;  (L_1, L_2,L_3) \;
 \eeqa 
 arose from the  variation of  (\ref{action}) with regard to infinitesimal changes of the point coordinates $\delta x_i$ and to infinitesimal rotations $\delta j_i$ of the {\em tri\`edres}.   This accorded to the venerated principle of  virtual velocities:
 \beqa \delta \int_{A_o} W dx &=& \int_{S_o}\sum_i (F_i \delta x_i + J_i  \delta j_i) \,d\sigma \\
  &+& \quad  \int_{A_o}\sum_i (X_i \delta x_i + L_i  \delta j_i)\, dx \; ,
 \eeqa
 where $d\sigma$ denotes the surface element on $S_o$ \citep[p. 597]{Cosserat:1909Note}.  By dissecting the medium  along a surface $S$ inside $A$ analog expressions could be discerned for the surface densities of {\em internal force }    $(F_i)$ and {\em torque} $(J_i)$  (``effort et moment de d\'eformation'').\footnote{The original notation was $(F,G,H)$ for our $(F_i)$, $I,J,K$ for our $(J_i)$, and $(X,Y,Z)$ for $(X_i)$, respectively $(L,M,N)$ for the $(L_i)$.}

 The evaluation of the invariance under Euclidean motions became  a complicated task. After diverse  transformations the authors derived two sets of equations involving auxiliary quantities   $p_{ij}$ and $q_{ij}$  (``dix-huit nouvelles auxiliaires'') which in slightly streamlined notation read as
 \beqarr X_j &=&   \sum_j \frac{\partial p_{ij}}{\partial x_i}  \label{equilibrium Cosserat 1}\\
 L_j &=&   \sum_j \frac{\partial q_{ij}}{\partial x_i} + p_{kl}-p_{lk} \label{equilibrium Cosserat 2} 
 \eeqarr 
where the $(j,k,l)$  in the last line form cyclical  permutions of $(1,2,3)$. Moreover they found that the internal forces and torques can be expressed by the auxiliary quantities:
\beq F_j = \sum_j p_{ij} n_i\; \qquad J_j = \sum_i q_{ij} n_i \, ; \label{internal force and torque Cosserat}
\eeq
here the  $n_1,n_2,n_3$ are the componentes of an interior directed unit normal of a surface element (of unit area) at any point
  \citep[p. 601, eqs. (29), (30)]{Cosserat:1909Note}.\footnote{Again our notation is slightly streamlined; the original notation was $p_{xx}, p_{xy}, \ldots p_{zz}$ for the stress densities  $q_{xx}, q_{xy}, \ldots q_{zz}$ for torque \citep[p. 601]{Cosserat:1909Note}. In the book \citep{Cosserat:1909Theorie} the equations appear on p. 137 in exactly the same form.}

Later readers would  read the  $p_{ij}$ and $q_{ij}$  as tensors of  the surface density of {\em internal stress} and  {\em internal torque} respectively (sometimes also ``proper'' stress and torque;  ``effort de d\'eformation'' and ``moment de d\'eformation'' for the Cosserats). The equations (\ref{equilibrium Cosserat 1}), (\ref{equilibrium Cosserat 2}) and  (\ref{internal force and torque Cosserat}) are now known as the {\em fundamental equations of elastostatics}.
Cosserats' theory includes the older linear theory of elasticity as a special case:  If the densities of external force,  torque and  internal momenta vanish, equation (\ref{equilibrium Cosserat 2}) implies $X_j = 0, L_j=0, q_{ij}=0$. Then   the stress tensor is symmetric, $p_{ij}=p_{ji}$, and satisfies  $  \sum_j \frac{\partial p_{ij}}{\partial x_i} =0$ which, in a dynamical context, may be  considered as a ``conservation condition''. Here, in the statical context, they indicate the equilibrium of the integrated forces acting on a closed surface inside the medium.

After having derived the equilibrium conditions for the surface torque momenta (see below), the Cosserats stated just that: 
\begin{quote}
Les auxiliaires que nous venons d'introduire et les \'equations qui les lient ne paraissent pas avoir \'et\'e jusqu'ici envisag\'ees sous une forme aussi g\'en\'erale; \`a notre connaissance, elles n'ont \'et\'e consider\'ees que dans les cas particuli\`er o\`u les neuf quantit\'ees $q_{xx}, \ldots,q_{zz}$ (Cosserat's expression for the components of the surface density of torque, here denoted by $q_{ij}$, ES) sont nulles, et le premier travail qui traite alors de la question semblable \`etre celui de M. Voigt  \citep[p. 137]{Cosserat:1909Theorie}.\footnote{``The auxiliary functions that we just introduced and the equations that relate them do not appear to have been envisioned in a form that was that general up till now; to our
knowledge, they have been considered only in the particular case in which the nine quantities $q_{xx} ,\ldots , q_{zz}$  are null, and the first work to treat that question seems to be that of {\em Voigt} \citep[p. 132,  Delphenich's English translation]{Cosserat:1909Theorie} . Moreover they recommend to compare with the papers  \citep{Larmor:1891,Love:1892/1906,Combebiac:1902}. In part of the literature P. Duhem is mentioned as a possible source for the Cosserats' turning towards  oriented elements of the contiuum. This seems implausible, however, because they did not mention him at this place, but only later with regard to the use of reversible transformations \citep[p. 73f]{Cosserat:1915}. The other way round, Duhem quotes the Cosserats positively in his \citep[p. 3]{Duhem:1906}; see \citep[pp. xxv, xxxv]{Brocato/Chatzis:Cosserat}. \label{fn Voigt-Cosserat}} 
\end{quote}
In a footnote added by the Cosserats they cited \citep{Voigt:1887,Voigt:1900}.
 
Given the generality of the assumptions, this was a great achievement. But the great complexity of the calculations   made the results extremely difficult to absorb.\footnote{Readers interested in technical  details  may like to consult  \citep{Badur/Stumpf}.} Only in retrospect could the Cosserats' theory be put into the context of wider mathematical theories and their derivations be justified  on the basis of general theorems, which involved less, or at least different, calculations: The equations (\ref{equilibrium Cosserat 1}) and (\ref{equilibrium Cosserat 2}) were identified as the Noether equations with regard to translational, respectively rotational invariance of the action \citep{Hehl/Obukhov:2007}. \citet{Pommaret:1997} sees them as special case of non-linear Spencer transformations in the theory of partial differential equations. Elasticity theorists had to develop their own viewpoint which gave reasons to address the  study of general elastic media.  In any case, 
Cosserat theory did not enter the broader theoretical or even experimental research for at least half a century. It was revived only in the 1950/60s. Even today it cannot  be considered mainstream, although it now seems to form  an  interesting sidestream of its own \citep[pp. xi ff.]{Brocato/Chatzis:Cosserat}.

During the course of their work, the Cosserats developed a perspective of  a grand unifying scheme for theoretical mechanics, covering  hydrodynamics,   heat conduction, electrodynamics, and elasticity \citep{Cosserat:1915}.\footnote{For  a detailed discussion of this point see \citep[sec. 4, pp. xxxvi ff.]{Brocato/Chatzis:Cosserat}}
This turned out to be  an  untimely enterprise: the relativity theories, special and general, and the rising quantum mechanics were just changing the role of  rational mechanics in mathematical physics. Although classical mechanics was not  invalidated in its core, it lost its central and foundational role for  natural philosophy of the 20th century. In consequence the overarching perspective of  Cosserats'  research program lost much of its power of persuasion. This may have contributed to the relative neglect of their generalized elasticity theory, in addition to its intrinsic technical difficulties.

On the other hand, the  theory of elasticity proposed by E. and F. Cosserat was highly valued by a small group of mathematical scientists, mainly in France but also internationally.\footnote{\citep[sec. 5, pp. xxxivff.]{Brocato/Chatzis:Cosserat}}
 Our protagonist, \'Elie Cartan, was one of the admirers. His re-reading of the  Cosserats' elasticity theory took place in the wider context of his investigations of Einstein's gravity theory,  which we consider next.

\section{\small Cartan's re-reading of Einstein gravity \label{section Cartan-Einstein}}
As already remarked,  Cartan's new geometric ideas were spelled out at the occasion of his studies of Einstein's general theory of relativity.  He started with analyzing  the form of the Einstein equation from a mathematician's point of view.  In  Riemannian geometry,  with metric $g= \sum_{ij}g_{ij}d_idx_j$, it may be written summarily as
\beq G  = \kappa T\; , \label{Einstein equation}
\eeq 
where $T$ denotes the energy-momentum-stress tensor of matter, $\kappa$ the gravitational constant ($\kappa = 1$ for Cartan).\footnote{Following Einstein, Cartan used a different sign $G=- T$. This is a question of conventions,  expressing the choice of a different sign  for the Ricci contraction of the Riemann tensor and the signature dependence of energy momentum. } 
\begin{itemize}
\item[(o)]   $G= \sum_{ij} G_{ij}dx_idx_j$, abbreviated $G_{ij}$, is a symmetric covariant  2-tensor, which contains the first and second partial derivatives of the  components $g_{ij}$ only. Equivalently it can be percieved as a vector valued 1-form $G^i_{\;j}$.
\end{itemize}
In Einstein gravity the left hand side is the {\em Einstein tensor}, 
\beq G=Ric-\frac{R}{2}g \; , \label{Einstein tensor}
\eeq
 where  $Ric$ and $R$ stand for  the Ricci, respectively scalar curvature of the Levi-Civita connection associated to $g$.  
In Cartan's view, the study of gravitational equations (plural!) boils down to the question which covariants may serve   on the left hand side of eq. ( \ref{Einstein equation}) as the (non-linear) partial differential operator on $g$. In any case, one should take into account two constraints which  Einstein had emphasized as basic principles:
\begin{itemize}
\item[$\;$(i)] $G$ is linear in the second partial derivatives $\partial \partial g$,
\item[(ii)] $G$ satisfies the {\em conservation law} (''loi de conservation'').
\end{itemize}
 (i) is necessary for avoiding too complicated differential equations.
 (ii) is a consequence of demanding a vanishing  covariant divergence of the energy momentum tensor. By  (\ref{Einstein equation}) this translates to the left hand side as $ \nabla_i G^i_j = 0$   in Ricci calculus (with $\nabla$ the covariant derivative associated to $g$).  Cartan preferred to  express conservation as the vanishing of  exterior covariant differential of $G$, which we denote hear as 
\beq d_{\omega} G = 0 \;,
\eeq  
because  it is defined with regard to a Cartan connection $\widetilde{\omega}=(\omega^i, \omega^j_{\;k})$.\footnote{For the torsion-free case see \citep[199]{Cartan:1922a}; for the general case \citep{Cartan:1923/24}; modern presentations in, e.g., \citep{Hehl/McCrea:1986}, \citep[pp. 269ff.]{Gasperini:2017} etc.} 

In his first paper on Einstein gravity Cartan  introduced his method of  differential forms   (outlined in section 1) for Riemannian geometry only. Using the Cartanized  coefficients  of the Riemann curvature (eq. (\ref{equation courbature})) he showed that  a  $G$ satisfying conditions (o) and (i) is  a  linear combination of  $Ric$, $ R \, g$ and $g$, all three expressed in terms of the basic differential forms $\omega^1, \ldots \omega^n$ \citep[p. 196]{Cartan:1922a}. If also the constraint (ii) is taken into account only the Einstein tensor form (\ref{Einstein tensor}) plus a linear  term  $ g$ remains, in the symbols introduced above:
 \beq G = \alpha (Ric- \frac{R}{2} g) + \beta  g \; , \label{general form Einstein tensor}
 \eeq
with two arbitrary constants $\alpha, \beta$ \citep[p. 203]{Cartan:1922a}.

From a mathematician's point of view, that was a highly pleasing result. Cartan was cautious, however, whether something similar had not perhaps been  already derived (in terms of the Ricci calculus) and published elsewhere in the international literature which, due to the effects of the great war, may have remained unknown in Paris.\footnote{``\' Etant donn\'ee la difficult\'e qu'on rencontre \`a avoir connaissance des M\'emoires parus \`l'\'etranger pendat la guerre et depuis la guerre, je ne suis pas absolument s\^ur qu'aucune d\'emonstration de ce th\'eor\`eme n'ait \'et\'e donn\'ee'' \citep[p. 142]{Cartan:1922a}.\\
{\em Taking into account the difficulty for gaining knowledge of foreign publications  during or after the war, I am not absolutely sure that no demonstration of this theorem has perhaps already been given.}
 } 
In fact, more or less at the same time at wich  Cartan wrote his manuscript of \citep{Cartan:1922a} H. Weyl proved that in Riemannian geometry the scalar curvature $R$ is the only invariant  containing not more than the first and second derivatives in $g$, and the second ones  only linearly. The proof was published about the time of Cartan's submissions of his notes to the {\em Comptes Rendus}  in an  appendix to the fourth edition of  {\em Raum - Zeit - Materie}  \citep[p. 287f., Anhang II]{Weyl:RZM4}.\footnote{In the French translation \citep[p. 279f.]{Weyl:RZMfranz}.} In the framework of a Lagrangian approach  Weyl's theorem  implied the same restriction for the Einstein tensor, which Cartan had derived.\footnote{Condition (ii) is here a result of the contracted Bianchi identities.}
 Weyl's proof had the advantage of being much shorter, but Cartan's analysis went deeper to the basic principles and was  more general, independent of a Lagrangian approach to the Einstein equation. 

In our representation of the general form of  (\ref{general form Einstein tensor}) we have assimilated Cartan's result to the more common notation of tensor calculus. But we have to keep in mind that
  Cartan used a different mathematical representation. That influenced also his  interpretation of the Einstein tensor:
  \begin{quote}
  Nour regarderons ses composantes comme des coefficients entrant dans l'expression de la projection sur une direction fixe d'une tension appliqu\'ee \`a un \'el\'ement \`a trois dimension de l'univers \`a quatre dimensions  \citep[p. 199]{Cartan:1922a}.\footnote{Similarly, in intuitive description in \citep{Cartan:1922[57]}.\\
  {\em We consider its components as the coefficients  appearing in the expression of the projection along  a fixed direction, applied to  tension exercised on an element  of three dimension in the four-dimensional universe.}}
  \end{quote}
  In other words, he conceptualized  the ``tenseur gravitationel''   $G$ as a {\em vector valued}   (alternating) {\em 3-form} $\widetilde{G}$, which, by analogy to classical elasticity, expresses the respective  stress force (``tension'') exercised on a 3-dimensional volume element in the 4-manifold (``l'univers'').  From a later point of view $\widetilde{G}$ may be understood in the Riemann geometric view as the Hodge dual of $G$.
  
  Following Cartan we shall use the terminology {\em  gravitational tensor} 
  (``tenseur gravitationel'') or {\em Einstein form} and the notation $\widetilde{G}$ if we conceptualize it  as a $(n-1)$-form (in dimension $n$), while the {\em Einstein tensor} (notation $G$)  will generally be understood as a the symmetric covariant tensor with coefficients $G_{ij}$.

To understand better what Cartan meant, we have to go  into more detail. Cartan decomposed $\widetilde{G}$ into its (real-valued) 3-form components $\Pi_i$,  such that it may be written as
\beq \widetilde{G}= \sum_i e_i \Pi^i \;  \label{Einstein form 0}
\eeq
($e_i$ the basis vectors of the Cartan orthonormal frame).\footnote{\citep[p. 203]{Cartan:1922a}, \citep[(vol 41) p.13]{Cartan:1923/24}} 

 Close to the end of analyzing the general form  of $\widetilde{G}$ he defolded  an intriguing argument involving his representation of the Riemannian curvature in terms of rotation coefficients $A^i_{\;jkl}$ (eq. (\ref{equation courbature})) and found that the  components $\Pi^i$ can be written as:\footnote{Up to the  factor $\alpha$ in eq. (\ref{general form Einstein tensor}) and  an equivalent to the ``cosmological'' term $\beta g$ which we suppress here.} 
  \beq \Pi^i = \sum \epsilon_k\epsilon_l\,  sgn(i,j,k,l)\, [\omega_j\,\Omega_{k l}] \;  \label{Einstein form 1} \\ 
  \eeq
where the system of indices $(i,j,k,l)$ is any cyclic permutation of $(0,1,2,3)$,  $sgn(i,j,k,l)$ its sign, and the summation runs over all such cyclic permutations \citep[p. 203]{Cartan:1922a}.  This form of the gravitational tensor  reappears in \cite[(vol. 41) p.13f.]{Cartan:1923/24}, where he called it the  {\em kinetic quantity of the mass} ``quantit\'e de mouvement masse''.\footnote{Cf. A. Trautman's commentary in \citep[p. 17]{Cartan:1986}.} Adopting the signature convention $sig\, g = (+---)$ for the metric, Cartan gave the ``kinetic quantity of mass'' (the  gravitational tensor)  in complete form, which displays its character as a vector valued 3-form openly  (ibid. eq. (7')):
\beq \widetilde{G}= \sum_{(ijkl)} sig(ijkl)\,  e_i\,[\omega_j\Omega_{kl}+ \omega_k\Omega_{lj}+ \omega_l\Omega_{jk}]\;   \label{Einstein form 2}
\eeq
(with summation over all cyclic permutations $(ijkl)$ of $(0\ldots 3)$).\footnote{In this  formula  Cartan wrote $[me_i]$ in the place of $e_i$, apparently to make the point dependence of the basis vectors $e_i$ immediately visible in his notation. \label{me-i} }

During the course of his study of the  gravitational tensor Cartan started to think geometrically about it and, as a consequence of the Einstein equation, also about the matter tensor. In the note of February 13, 1922, he stated:
\begin{quote}
On sait que, dans la th\'eorie de la relativit\'e g\'en\'eralis\'ee d'Einstein, le tenseur qui caract\'erise compl\`etement l'\'etat de la mati\`ere au voisinage n'un point d'Univers et {\em identifi\'e} \`a un tenseur faisant intervenir uniquement les propri\'et\'es {\em g\'eom\'etriques} de l'Univers au voisinage de ce point \citep[p. 437, 1-st emphasis ES, 2-nd emphasis in the original]{Cartan:1922[57]}.\footnote{\em One knows that in Einstein's generalized theory of relativity the tensor which  completely  characterizes the state of matter in a neighbourhood of the Universe is \underline{identified} with a tensor which is exclusively made up by the \underline{geometric} properties of the Universe at this point. ({\em 1-st emphasis ES, 2-nd emphasis in the original}) } 
\end{quote}
This differed from how  physicists usually understand the Einstein equation. For them  eq. (\ref{Einstein equation})  expresses a kind of communication between two aspects of reality, spactime and matter, not a reduction of one to the other. Einstein fought strongly against the claim that his theory of general relativity had geometrized gravity.\footnote{\citep{Lehmkuhl:2014}} 
But for Cartan the idea that the Einstein equation justifies an identification of  its  left hand (geometrical) side and the right hand (matter) side  became  a guiding motif for his further investigations.  In his notes of 1922 he used the notions  ``tenseur de mati\`ere'' ``tenseur d'energie d'Einstein'' etc.  synonymously and understood them to be {\em defined} by  geometrical curvature properties.\footnote{This identification is announced already in  the title of \citep{Cartan:1922[57]} and referred to in the next notes, e.g.  \citep[593]{Cartan:1922[58]}.} 
This was a clue for his way of  generalizing  Einstein gravity in \citep{Cartan:1923/24,Cartan:1925} and, to my knowledge, remained so in the years to come.

In order to explain what he meant Cartan used a 3-dimensional analogue of the Einstein equation. Then the right hand side reduces to the classical (symmetric) stress tensor of matter, and the left hand side analogue may be described by curvature properties of a space with the correct properties of infinitesimal  tri\`edres (3-frames) which, if one wants so, define their own metric different from the classical Euclidean metric of the ordinary embedding space. In Cartan's conceptualization of curvature the latter expresses itself in a  ``rotation compl\'ementaire'' (complementary rotation) which is to be  applied after parallel transporting  a tri\`edre  around an infinitesimal loop. Presupposing the quid pro quo mentioned in the our introduction as self-evident he continued:
\begin{quote}
Cette rotation peut se representer par un vecteur. L'\'etat de divergence entre l'espace donn\'ee et l'espace euclidien peut donc \^etre traduit par un vecteur attach\'e \`a chaque \'el\'ement de surface orient\'e de l'espace. \citep[p. 438]{Cartan:1922[57]}\footnote{\em This rotation may be represented by a vector. The state of divergence between the given space and Euclidean space can thus be expressed by a vector attached to any oriented  surface element of the space.}
\end{quote} 
In the same note he declared that the assignment of vectors to surface elements  results in a tensor from which one can show symmetry and ``conservation law'' just like for the original Einstein equation.  In 3-dimensional (Euclidean embedding) space the expression of an infinitesimal rotation by a vector was a standard procedure. 
Cartan concluded:
\begin{quote}
Il r\'esulte de ce qui pr\'ec\`ede qu'on peut expliquer l'\'etat d'un milieu \'elastique en \'equilibre en admettant que l'espace qui le contient est d\'eform\'e et qu l'\'etat de tension du milieu  traduit physiquement cette d\'eformation g\'eom\'etrique. (ibid.)\footnote{\em From the preceding it follows that one can explain/express the state of an elastic medium in equilibrium by assuming that the space in which it is contained is deformed and that the state of tension of the medium reflects this deformation physically.  }
\end{quote}
This is an interesting sentence. For the moment  we leave it open  whether we ought to understand ``expliquer'' in the sense of making something explicit in a mathematical sense, or even stronger as an explanation in the physical sense. 

In order to understand what the mathematics behind this sentence is, one has to see the context. Immediately after this discussion of 3-dimensional classical elasticity, Cartan explained a geometrical interpretation of the Einstein equation (\ref{Einstein equation})
in the light of (\ref{Einstein form 1}) derived in \citep{Cartan:1922a}.\footnote{Remember that the \citep{Cartan:1922a} was already written at the time of submission of the {\em Comptes Rendus} notes.} Although he did not discuss  elasticity in the latter, his notes  show that  in early 1922  he thought about {\em classical elasticity as a 3-dimensional analogue of the Einstein equation}. In this case  the right hand side reduces to the tension tensor which Cartan would understand as a vector valued 2-form with real valued 1-forms  $\tilde{T}^i$ as components ($i=1,2,3$). It expresses  the stress force $\tilde{T}^i(\sigma)$ exercised on any infinitesimal surface element $\sigma$.

In a 3-dimensional version of  (\ref{Einstein form 0}), (\ref{Einstein form 1}) the gravitational 3-form reduces to a 2-form and the signature coefficients are all $\epsilon_j=1$.\footnote{If we consider the Einstein tensor in dimension $n=3$ and allow us (anachronistically) to apply Hodge duality, we  arrive at $\widetilde{G}$ as a {\em vector valued 2-form}.}
 With Cartan's choice $\kappa = 1$ a 3-dimensional analogue of the Einstein equation would be
  \beq \sum_{k,l}  sgn(i,k,l)\, \widetilde{\Omega}_{k l} = \tilde{T}^i \; \label{3-dimensional Einstein equ}
    \eeq
with $ \widetilde{\Omega}_{k l}$ the componentes of the curvature 2-form of a 3-dimensional Cartanized Riemannian geometry. 
The alternating signs in the  summation of (\ref{3-dimensional Einstein equ}) associate a vector to the  rotational coefficients just like in the vector product representation of infinitesimal rotations.\footnote{For $\tilde{\Omega}= \left( \begin{array}{ccc}
0 & c & -b \\ 
-c &0  & a \\ 
b & -a & 0
\end{array}\right)  $ eq. (\ref{3-dimensional Einstein equ}) gives $T =\left( \begin{array}{c}
a \\ b \\ c \end{array} \right) $.} 
This would underpin what Cartan  intuitively  circumscribed in his  note \citep{Cartan:1922[57]} quoted above and  explain half of the {\em quid pro quo} cited in our  introduction (p. \pageref{quid pro quo}), the expression of a rotation (curvature) by a vector.

Two years later, in the second lot of \citep{Cartan:1923/24}, Cartan gave a more technical explanation in terms of  the 4-dimensional Einstein equation. Here he   concentrated  on a 3-dimensional  spacelike hypersurface $S $ (in the infinitesimal  represented by a hyper{\em plane}) corresponding to $x_o=0$ and ``projected''  the 4-dimensional rotations onto $S$.  If a vector $\xi = (\xi^i)$ is rotated by  $\Omega^i_j$, the hyperplane projection of the change  is $\Delta \xi^i = \Omega^i_j \xi^j$. The rotation of the hyperplane itself is therefore given by the components $(\Omega_{23}, \Omega_{31},\Omega_{12})$. Cartan continued:
\begin{quote}
Elle peut, dans ce hyperplan, \^etre repr\'esent\'ee par le bivecteur \\[-0.8em]
\[(\ast)  \qquad [e_2e_3]\Omega^{23} + [e_3e_1]\Omega^{31}+[e_1e_2]\Omega^{12}\]
 our encore par le vecteur polaire de m\^eme  mesure \\[-0.8em]
 \[ (\ast \ast) \qquad \frac{1}{\sqrt{g_{11}g_{22}g_{33}} } (e_1 \Omega_{23} + e_2 \Omega_{31} + e_3 \Omega_{12}) \, \]
  \citep[(vol. 41) p. 16, (marks $(\ast), (\ast\ast)$ added, ES)]{Cartan:1923/24}.
\end{quote}

Here  ($\ast$) was a Grassmann type characterization of (infinitesimal) rotations, while ($\ast\ast$) was its equivalent in terms of a vector product representation.  This   transition  from ``bivectors'' to ``polar vectors''  was a special case of what  Grassmann had introduced as a more general duality (called ``Erg\"anzung'' by him).\footnote{In modernized notation, Grassmann established an equivalence (isomorphism)  between $\Lambda^k V$ and $\Lambda^{(n-k)}V$ for any $n$ dimensional {\em Ausdehnungsgebiet} ($0 \leq k \leq n$) with a  volume form, respectively  a basis $e_1,\ldots,e_n$ with the property $e_1\wedge\ldots \wedge e_n =1$ (using modern notation for alternating products).  He assigned to basis elements $e_{i_1}\wedge \ldots \wedge e_{i_k}$ ($i_1 < \ldots < i_k$) in $\Lambda^{k}V$ the bais elements $e_{j_1} \wedge \ldots e_{j_{n-k}}$ ($j_1 < \ldots < j_{n-k}$) of $\Lambda^{(n-k)}V$ for which $e_{i_1} \wedge \ldots \wedge e_{i_k}\wedge e_{j_1}\wedge \ldots \wedge e_{j_{n-k}}= e_1 \wedge \ldots \wedge e_n =1$ and used linear continuation 
\citep[\S 89, p. 57f.]{Grassmann:1862}, cf. \citep{Scholz:1984Grassmann}.  As Grassman introduced an inner product in $V$, which made $e_1, \ldots, e_n$ an orthonormal basis, this can be considered as a linear algebraic isomorphism serving as the basis for the  later {\em Hodge duality}.}
In dimension 3  this type of dualization was particularly  well known, even without any reference to Grassmann. 
But Cartan was well aware of the  general nature of   Grassmann dualization and  used it in his discussion of the  invariants of his geometry \citep[vol. (40) p. 400ff.]{Cartan:1923/24}.

From this point if view the first half of Cartan's quid pro quo 
resulted from a {\em Grassmann type duality transformation} (of infinitesimal rotations in dimension $n=3$ to polar vectors), which Cartan introduced for gaining a {\em geometrical understanding} of the spacelike part of the  {\em Einstein equation}. The result of the transformation led to a new geometrical picture of {\em classical elasticity}: stress could be expressed in terms of the curvature of a space with a connection and metric adapted to the mechanical properties of the material medium under investigation.

\section{\small   Einstein gravity in analogy to  geometrized Cosserat  elasticity  \label{section Cosserat-Einstein} }

The second part of the {\em quid pro quo} resulted from Cartan's generalization of Einstein's theory of gravity and involved an adaptation of Cosserat elasticity to his research program of 1921/22.  At the time of submitting his {\em Comptes Rendus} notes, in February and March 1922, Cartan had all this  in mind, but it took some time to work out the mathematical details. They are contained in the two-part paper \citep{Cartan:1923/24,Cartan:1925} the first part of which came in two lots (vol. 40, 41 of the {\em Annales ENS}).\footnote{Reprint in \citep{Cartan:1955}, English translation with a  commentary (foreword) by A. Trautman in \citep{Cartan:1986}.} 

 In this paper he showed that the vacuum Maxwell equations are compatible with any (Cartan-)  connection of the Poincar\'e group; but taking Lorentz forces into account may run into difficulties. For a kinetic quantity of energy  like in (\ref{Einstein form 2}) the Lorentz  force exercised on an electric current density  came out correctly, i.e., in agreement with special relativity,   only if the ``universe'' has vanishing torsion \citep[p.20f.]{Cartan:1923/24}. 
That was  disappointing, but Cartan indicated a way out:
 \begin{quote}
 La conclusion pr\'ec\'edent (vanishing torsion, ES) ne serait pas logiquement n\'ecessaire si l'on admettait une conception de la M\'ecanique des milieux continus plus large que la conception habituelle, la ``quantit\'ee de mouvement-masse'' \'el\'ementaire \'etant rep\'esent\'ee par une syst\`eme de vecteur et de bivecteurs\\[-0.8em]
 \[ G = [me_i]\, \Pi^i + [e_i e_j]\, \Pi^{ij} \, 
 \] 
 \citep[p. 21]{Cartan:1923/24}.\footnote{``The above conclusion is not logically forced upon us if we accept a broader framework for mechanics of continuous media and
represent the energy-momentum density by a system of vectors and
bivectors: $ G = [me_i]\, \Pi^i + [e_i e_j]\, \Pi^{ij} $ '' \citep[p. 123]{Cartan:1986}.
 
 Compare fn. \ref{me-i}.
 }
 \end{quote}
 
If such a modification of the Einstein form is accepted,  the  laws of electromagnetism, including the Lorentz forces,  were compatible with a non-vanishing torsion. In this context, it was natural to assume that the ``quantit\'e de mouvement-masse elementaire'', i.e. $\widetilde{G}$, should remain a     geometric integral invariant, like in Einstein's theory.
This, so Cartan declared, was easy to achieve. One had only to replace 
 the rotation associated to any surface element by the total displacement of the full Cartan curvature (``d\'eplacement total ({\em rotation et translation})'') assigned to the surface element. 
 
The rotations had been transmuted to vectors  by Grassmann duality in 3-dimensional spacelike hyperplanes and this transmutation was taken over to $n=4$ (see above). In an analogous manner Cartan transmuted  translations into bivectors (Grassmann duality in the 3-dimensional spacelike projection, but here  transferred to  $n=4$). 
  In this way  $\widetilde{G}$  became a  fully Cartanized variant of the Einstein form  \citep[p. 22, eq. (11)]{Cartan:1923/24}:
 \beq \widetilde{G}= \sum_{(ijkl)} sig(ijkl)\,\left(  e_i\,[\omega_j\Omega_{kl}+ \omega_k\Omega_{lj}+ \omega_l\Omega_{jk}] - [e_i e_j] [\omega_k\Omega_{l}-\omega_l\Omega_k] \right) \label{Einstein form 3}
 \eeq
 
Expressed in more recent terminology Cartan proposed a  3-form with values in the Grassmann algebra of the tangent bundle as a generalization of the gravitational tensor. It  consists of two terms, the first one with values in $TM$\footnote{In fact, Cartan put square brackets about his symbols of the vector basis $[me_i]$, apparently in order to emphasize the Grassmann character of the term.}
 contracts rotational curvature and transmutes it into a vector.  
 The second one with values in the bivector bundle  $\Lambda^2 (TM)$ {\em transmutes translational curvature  into a bivector}.
 
 That seemed to agree nicely with the structure of Cosserat elasticity. After  considering his geometrical interpretation of the Einstein form as a new representation  of classical elasticity of a (1-parameter) homogeneous medium, Cartan was  now  tempted to read the 3-dimensional reduction of his generalized Einstein equation as a geometrization of Cosserat elasticity. The vector part of $\tilde{G}$ resulted, after Grassmann dualization, from the rotational component of Cartan curvature and could express the stress forces on surface elements. Similarly due to Grassmann dualization, the bivector component  resulted from the translational curvature and could be be interpreted as a rotational momentum exercised on surface elements. Cartan even gave an argument  that, due to the specifique form of the torsion in his case, the translation associated to any surface element is normal to the latter \citep[p. 21, fn]{Cartan:1923/24}.\footnote{At another place he showed that this property implies the identity of autoparallels of the Cartan connection with the geodescis of the related Riemannian structure \citep[\S 66, p. 407]{Cartan:1923/24}.}
 This was apparently the technical insight standing behind Cartan's verbal description of Cosserat elasticity in his note \citep{Cartan:1922[58]}.  Moroever, it contains a resolution of the riddle of the  quid pro quo. But here a warning is appropriate:  Cartan no longer mentioned the Cosserat theory in his  \citep{Cartan:1923/24}. We come back to this point in the next  section.

 Cartan encountered another problem which appeared more serious to him. 
In general, the covariant exterior derivative of the generalized Einstein form does not vanish.  Our author considered this as an important criterion for a medium in equilibrium and proposed to ensure it by 
 imposing some appropriate  constraint, the easiest one being  (ibid., p. 22)
\beq \sum (ijkl) \left( \Omega_j\Omega_{kl}+ \Omega_k\Omega_{lj} + \Omega_l\Omega_{jk} \right) = 0 \qquad \mbox{for all $i=0,\ldots,3$.} \label{Cartans constraint}
\eeq 
Later authors ( see section 7) replaced the condition of vanishing covariant divergence, the ``conservation law'' in Cartan's view, by the more general one of a contracted Bianchi identitiy for Riemann-Cartan geometry. Only for vanishing torsion or Cartan's algebraic constraint (\ref{Cartans constraint}) it boils down to the ``conservation  law'' \citep[p. 152f.]{Trautman:1973}.

Unhappily the relation (\ref{Cartans constraint}) gave an {\em algebraic} constraint for torsion  which   could therefore no longer play the role of a dynamical field in this approach. 
Perhaps this was a reason for  Cartan  to hesitate  claiming immediate  physical  relevance  for his theory.\footnote{Cf. the remark by A. Trautman in \citep[p. 17]{Cartan:1986}. \label{fn Trautman}} 
In any case, he emphasized:
 \begin{quote}
 On a ainsi une g\'en\'eralisation, au moins math\'ematique, de la th\'eorie d'Einstein, g\'en\'eralisation compatible avec toutes les lois de l'\'Electromagn\'etisme \citep[p. 22]{Cartan:1923/24}.\footnote{``This yields -- at least on a mathematical level -- a generalization of
Einstein's theory which is compatible with all the laws of electromagnetism'' \citep[p. 124]{Cartan:1986}.
 }
 \end{quote}
 
 He even indicated that this generalization of Einstein gravity can be formulated in a Lagrangian field approach. He proposed a generalization of the  Hilbert action expressed in terms of his geometrical quantities,
 \beq \mathfrak{L}_{grav} = \sum_{(ijkl)} (ijkl)\, [\omega_i\omega_j\Omega_{kl}] \; ,
 \eeq
but did not  consider a separate matter Lagrangian (in agreement with his ``unphysical'' identification of the left hand and the right hand side of the Einstein equation). 
 
 At this point, our author stopped his physics related studies without even shedding a first glance at the dynamical effects of his theory: 
\begin{quote}
\ldots  je me contente de cette indication, sans entrer dans plus de d\'etails \`a ce suject \citep[p. 23]{Cartan:1923/24}.\footnote{ {\em\ldots I  make to do with this indication without going into more details of this subject.}\\
Translation in \citep[p. 124]{Cartan:1986}: ``However, we shall not discuss this issue in more detail.'' }
\end{quote}
He rather continued with studies of  Weyl's dilational gauge metric and in part II,   according to the title of the paper,  with  studies of Cartan spaces of the affine group or subgroups of it \citep{Cartan:1925}. 

In other papers, e.g. \citep{Cartan:1924[70]},  he  turned towards generalizations of the group. Where he sticked to the Euclidean group he tended to investigate more classical questions. In particular he  showed that Clifford parallelism in elliptic space can be  seen, from the point of view of his new methods, as due to a  connection with torsion but vanishing rotational curvature, while the autoparallels coincide with elliptic straight lines \citep{Cartan:1924[71]}.\footnote{Cf. \citep[p. 46f.]{Cogliati/Mastrolia}.} All these topics turned his and his readers' attention away from the Cosserat inspired view of torsion. It rather  hinted in the direction of what a little later became the study of distant parallelism or {\em absolute parallelism}, as Cartan would call it.

\section{\small Cartan's practice of mathematical analogies  \label{section XX2}}

It is not completely clear,  why Cartan stopped short of pursuing the physical considerations    further. One obvious reason might have been  that he was more interested in studying the mathematical side of his approach than going deeper into  the  physics. But he seems also to have developed doubts  with regard to the physical relevance and the interpretation of his findings. The algebraic constraint for torsion which he perceived as necessary to satisfy the ``conservation law'' as he understood it must have been a stumbling block for him.\footnote{Cf. fn \ref{fn Trautman}.} 

We even find some hints in  \citep{Cartan:1923/24} that during the continuation of his work  he  developed  doubts regarding the feasibility of his  geometric interpretation  of  Cosserat elasticity in the {\em Compte Rendus} notes. 
 In \S 60 of his paper he  discussed the  differential forms  of  a (Cartan) space with the Euclidean group, invariant under changes of the Cartan reference system (change of Cartan gauge in modern terminology). One of the invariants combines a ``syst\`eme de vecteur et de couples'' of a  form similar to (\ref{Einstein form 3}):
 \beq  [m e_1]\Omega_{23}+ [m e_2]\Omega_{31} + [m e_3]\Omega_{12} + [e_1 e_3]\Omega_1 + [e_3 e_1]\Omega_2 + [e_1 e_2]\Omega_3\,. \label{analogie trompeuse}
 \eeq 
 This is a 3-dimensional analogue of the expression for the generalized Einstein form (\ref{Einstein form 3}). But here the relation between translations and bivectors is simpler and more direct than in the 4-dimensional case (in fact it is  nothing but the Grassmann duality between vectors and bivectors which only holds for dimension $n=3$). 
   Cartan commented:
 \begin{quote}
 Le syst\`eme  de vecteurs cet de couples 
 \[ [m e_1]\Omega_{23}+ [m e_2]\Omega_{31} + [m e_3]\Omega_{12} + [e_1 e_3]\Omega_1 + [e_3 e_1]\Omega_2 + [e_1 e_2]\Omega_3\,\]
repr\'esente de m\^{e}me  le d\'eplacement associ\' e \`a un \' el\'ement de surface. Sa d\'eriv\'ee ext\'erieure est nulle.  Si l'on regardait ce syst\`eme comme repr\'esentent les tensions qui s'ecercent sur un milieu mat\'eriel (tensions comportant des couples), ce milieu serait en \'equilibre \citep[vol. 40, p. 401, fn]{Cartan:1923/24}.\footnote{``The system
 $ [m e_1]\Omega_{23}+ [m e_2]\Omega_{31} + [m e_3]\Omega_{12} + [e_1 e_3]\Omega_1 + [e_3 e_1]\Omega_2 + [e_1 e_2]\Omega_3\,$
of vectors and torques represents the same displacement associated with a surface element.
Its exterior derivative is zero. If this system is viewed as describing tensions (with torque)
in a material medium, the medium would be in equilibrium.'' \citep[p. 96, fn. (9)]{Cartan:1986}
 }
 \end{quote}
 
 In other words,  an invariant of this type looked as if it could be interpreted as a geometrical expression for the system of tensions/momenta inside a Cosserat medium. But, different from his announcement in \citep{Cartan:1922[58]} of February 1922, in which he had  referred so positively to the 
 ``beaux travaux de MM. E. et F. Cosserat'',  Cartan now continued  with a methodological reflection  which stepped  back from  an interpretation in this sense and  came even close to  a  disassociation:
\begin{quote}
Il y a l\`a une des nombreuse analogies, plus ou moin trompeuses, qui existent  entre la G\'eom\'etrie et la M\'ecanique. En fait, ce n'est qu'une analogie. (ibid.)\footnote{{\em We have here  one of the numerous, more or less misleading, analogies  which exist between geometry and mechanics. In fact this is not more than an analogy.}
\\ The translation in \citep[p. 96, fn. (9)]{Cartan:1986}  (``We have here one of the numerous somewhat misleading analogies between geometry and mechanics.'' ) omits the last phrase of Cartan's remark.
}
\end{quote}

Cartan  neither claimed that the analogy {\em was} ``trompeuse'' (erroneous), nor did he repeat a positive claim of content for it. But it seems that he  was having second thoughts  about this point between February 1922 and the final preparation of the paper  \citep{Cartan:1923/24}. 
It now became clear to him that the analogy between Cosserat elasticity and the generalized Einstein equation was less perfect than he had inititally hoped  because of the non-vanishing covariant exterior derivative in dimension $n=4$. 
This would also explain why he avoided  any explicit reference to the Cosserats' work; he even did not  mention their name any longer in the new paper.\footnote{Later research on this question, starting in the 1960s, showed that a viable usage of Cartan geometry in Cosserat type theory of elastic media needed, in fact, a more sophisticated approach than was available to him in 1922 (see below).}

That could not change the crucial heuristic role which, according to Cartan's own testimony in the notes of early 1922, Cosserat elasticity  played   during  his early  work on the Einstein equation and its generalization in the light of his geometrical ideas.  Cosserat torque  would strongly underpin the transmutation from  translations to bivectors (rotational momenta) in 3-dimensions by a Grassmann type duality  (\ref{analogie trompeuse}), {\em if it could be considered as physical}.    In 1921/22 Cartan was  apparently impressed by the analogy to the transmutation of the rotations to vector-like stress in his  re-reading of the Einstein tensor. The analogy helped Cartan to structure his argumentation in which he tried to build a  bridge between Einstein gravity and geometry. Once the bridge was built in the form of our eq. (\ref{Einstein form 3}), the reference to Cosserat media could be downgraded to the status of a mere analogy without losing too much. 

Only the {\em terminological residuum} of Cartan's early heuristics remained unaffected: the translational curvature baptized under the impression of the 3-dimensional analogy  with Cosserat media  when it was still fresh and strong, continued to be called {\em torsion} and remains so until today. 

\section{\small The aftermath  \label{section XX3}}

This is not the end of the story. We should not finish our's without having a short glance at the reception and further developments connected to Cartan's early papers on Einstein gravity and Cosserat theory.
The early idea of a potentially  intimate connection between Cosserat elasticity and Einstein gravity did not play a role in the reception for many decades to come, while   Cartan geometry attracted the attention of mathematicians by other reasons. J.A. Schouten got interested in Cartan's proposals of torsion in his general studies of connections. He  contacted Cartan in 1924. The two mathematicians communicated on  linear connections in Lie groups and found that left and right translations in any Lie group lead to distant parallelism structures, i.e. connections with torsion but vanishing curvature, with autoparallels which agree with geodesics of the Riemannian  metric on the group manifold induced by the Cartan-Killing form. They even were able to show that, with the exception of the 7-sphere, the Lie groups are the only Riemannian manifolds with this property  \citep{Cogliati/Mastrolia}.

Distant parallelism became a ``hot topic''  at the end of the the 1920s, when Einstein started to study it as the framework for one of his attempts to unify gravity with electromagnetism. After Cartan reminded Einstein that this approach could be well framed in his geometrical method and had been mentioned by him in talks with him in 1922 \citep[p. 4]{Cartan/Einstein}, Einstein hurried to give credit to Cartan and accepted that his study of gravity in terms of  {\em  distant parallelism} (also called {\em teleparallel} gravity) used a specific type of Cartan geometry.
 In this setting the deviation of flat space was encoded in the torsion part of curvature only, while the rotational curvature was set to zero.\footnote{\citep{Goenner:UFT}, cf.  \citep{Cartan1929}. }
Also R. Weitzenb\"ock  had studied flat linear connections (vanishing Riemannian curvature) with torsion in the course of his study of differential invariants, i.e. in a pure mathematics context  \citep[pp. 317ff.]{Weitzenboeck:1923}. He did not relate this to Einstein gravity at the time but was keen to get acknowledgement from Einstein and published a note on on the topic in 1929 \citep{Sauer:Fernpar}. Neither Cartan nor Cartan geometry was ever mentioned by him. 
 
  Although Cartan himself had given an example of a teleparallel Cartan space in his note \citep{Cartan:1922[58]}, the general outlook of this example was a far cry
 from his early idea of  interpreting torsion by  rotational momenta as an additional feature of the gravitational field. In a way it even was  opposite to his proposal for the physical interpretation of translational curvature. But even so, the studies of  distant parallelism in the gravity  context  demonstrated  the  openness  of general mathematical structures for different  physical interpretations. Even those which were designed with definite physical interpretations in mind, like Cartan geometry of the Euclidean and Poincar\'e group, did not carry the mark of their original interpretation with them  as some sort of inbuilt, although perhaps hidden finality.  
 
We also have to be aware that, by a constellation of  historical contingency,  the late 1920s was also  the time in which quantum physicists started to realize that the new complex (wave) fields  could carry an internal rotational momentum, called {\em  spin} \citep{Pauli:Spin,Dirac:1928_I/II}. But at this time no author had the idea that this new internal torque-like momentum  might give new support to  Cartan's idea of torsion.

This changed only much later, in the 1960s. An important contribution  for  renewing the interest in Cartan torsion among physicists arose in the wake of the work of D. Sciama and T. Kibble in gravity theory. Without knowing it, the two authors independently reinvented much of the Cartan geometric field structures by considering what physicists call the ``localization'' of the Poinar\'e group \citep{Sciama:1962,Kibble:1961}. They found that the spin of elementary particle fields might play a role for a generalized theory of Einstein gravity which was close to what Cartan had anticipated in his early papers. The close relationship of their theory to Cartan geometry was not clear to Sciama and Kibble; but it was   soon made explicit by other authors, at first by 
F. Hehl in his PhD dissertation \citep{Hehl:Diss} and independently  by A. Trautman \citep{Trautman:1973}. 

A group of authors joined and extended this research program.\footnote{Much information on this development is collected in the  reader  \cite{Blagojevic/Hehl} which contains very helpful commentaries. For systematic surveys see  \cite{Trautman:2006,Hehl:Dennis}. } 
They realized that Cartan geometry offered a tailor-made geometric framework for  infinitesimalizing  (``localizing'' in the language of physicists) energy-mometum and spin  currents known from Minkowski space and special relativity.
The Cartan geometry of this approach was modeled on the  Poincar\'e group and has both, rotational curvature and torsion, like in Cartan's work of the early 1920s.  If the gravitational Lagrangian was chosen as closely as possible to the Hilbert action of Einstein gravity, it turned out to be the one Cartan mentioned  in his side remark quoted above. Only Cartan's constraint (\ref{Cartans constraint}) had to be relaxed \citep{Trautman:1973}. The resulting theory has therefore correctly been called {\em Einstein-Cartan gravity}. It is  considered as a ``viable''  alternative to Einstein gravity   -- although one which can be distinguished from the latter only under the conditions of extremely high energy densities.\footnote{\citep{Trautman:2006,Hehl:Dennis}.} 

Motivated by the  successes of non-abelian gauge field theory in the rise of the {\em standard model} of elementary particle physics,   gauge field dynamics related to  the Poincar\'e group was studied in the 1970/80s. It turned out that Einstein-Cartan gravity can be reconstructed in a Cartan space modeled on the Poincar\'e group, where the dynamical equation for translational curvature --  paradoxically still called ``torsion'' --  couples to the energy-momentum current (of matter and the gravitational field), and the rotational curvature is coupled to the spin-current \citep[p. 337ff.]{Hehl_ea:1980Poincare-gauge-field-theory}. This was an intriguing result. Conceptually it set the couplings right while avoiding the surprising crossover of translational and rotational aspects for geometry and dynamics of Einstein gravity as seen by Cartan in analogy to Cosserat elasticity. In this way it generalized the ``teleparallel'' representation of Einstein gravity by adding spin.  But the authors presented the respective Lagrangian as a special case, even  a degenerate one, and directed their attention towards more general quadratic Lagrangians in Poincar\'e gauge field theory.\footnote{See also the surveys in \citep[p. 164f.]{Hehl:Dennis} and \citep{Blagojevic/Hehl}.}

This was not the only path which led back to Cartan's ideas of the early 1920s. Also Cartan's spare remarks on Cosserat elasticity found  successors a generation later, although only after a specific turn taken by authors interested in the mathematical  study of the yielding of plastic materials.
 K. Kondo in Japan  proposed to model dislocation in crystal matter by the torsion of a linear connection \citep{Kondo:1952}. 
   E. Kr\"oner, an expert in classical solid state physics, and F. Hehl's PhD advisor at the {\em Bergakademie/TH} Clausthal-Zellerfeld,  was attracted by the  ideas of the Cosserats on generalized elasticity and its link to Cartan's generalized geometry. 
   He was one of those authors who brought in the idea that such a material structure could be studied in a Cartan geometry modeled on the Euclidean group, called {\em Riemann-Cartan space} by these authors.\footnote{\citep{Kroener:1963torque,Kroener:1963}  }
This led to the attempt  of studying  dislocations and proper tensions in metals in terms of  Cartan-type geometrical methods and  was also the background from which Hehl  entered gravity.
 
    A complicated story started; we have to cut it short. 
      The different components of the Cartan connection  had to be linked to physical quantities expressing the deformation of the material, different types of dislocation inside the material, and the hypothetical force and torque stresses.\footnote{A recent survey of the resulting theory can be found in \citep{Hehl/Obukhov:2007}.} It turned out that a Cartan geometry with translational curvature but no  rotational curvature, thus with ``distant parallelism'' in the language of gravity theorists was an  option.  The intuitive idea behind this  was the choice of a Cartan reference system without rotational curvature, adapted to the geometry of the lattice structure of the material. Torsion can be expressed as  a {\em closing defect}, which arises if one    parallel  transports a vector $u$ along an infinitesimal shift vector $\delta x=v$ and  $v$ along $\delta' x= u$.\footnote{The  ``would-be'' infinitesimal parallelogram arising from this procedure does not close (thus it is {\em no} parallelogram). Mathematically, the defect is expressed by the asymmetry of the corresponding linear connection $\Gamma^i_{jk}$ in the lower indexes, $\Gamma^i_{jk}\neq \Gamma^i_{kj}$.} Such closing defects can be used to model dislocations in the material, once they turn up sufficiently often (densely). This results in a mathematical description of materials with densely distributed dislocations by a ``teleparallel'' Cartan structure. Its torsion, here clearly to be interpreted in its translational connotation,  is related to the dislocation field.\footnote{Cf \citep[p. 167ff.]{Hehl/Obukhov:2007}.  In this context, not any Riemann-Cartan geometry with distant parallelism is feasible. Additional constraints have to be observed, in order to  lead to acceptable deformation quantities related to the Cartan structure  \citep[p. 163]{Hehl/Obukhov:2007}.} 
      
     Finally the tables have been turned another time. At the origin of Cartan's theory the context of  continuum mechanics gave him the motivation for introducing his slightly paradoxical terminology of ``torsion''. Now we find an epistemic constellation even in the  field of continuum mechanics,  for which the geometrical naming of {\em translational curvature} would be closer to the matter than the terminology chosen by Cartan. But in the meantime the latter has been widely  established in the community.

In the end, Cartan's papers of the early 1920s have found new readers also among  present day  theorists of continuum mechanics. In this recent development it became clear that a reliable connection between the physics of matter and geometry needs much more sophistication than could be imagined by Cartan (or even the Cosserats). It turned out that, in the long run,  Cartan's analogy between  generalized elasticity and  gravity, which became apparent in the framwork of his geometry, was not as misleading (``trompeuse'') as he may have thought in 1924. The mathematical analogy established by Cartan   became  a stimulating input for these studies, even though the structural analogy  had to be disentangled, before it could bare fruits. 
\\[3em]

{\em Acknowledgments:} {\small This paper has been written in the honor of H. Sinaceur.  I thank the organizers of the workshop {\em 30 ans Hourya Sinaceur: ``Corps et mod\`eles''}, Paris, June 15--17, 2017, for the opportunity to present my findings on the  Cartan's early work  on gravity close to the place of his activity. The paper will appear in  French in the  proceedings of the workshop.  I thank Friedrich Hehl for  critical remarks and hints, in particular his  personal informations on the background of the Kr\"oner research group in the 1960s, Alberto Cogliati for making me aware of his research on the correspondence between Schouten and Cartan.  Emmylou Haffner generously took on the burden of translating my paper from  English to French.} 
\small
\newpage
\addcontentsline{toc}{section}{\protect\numberline{}Bibliography}

 \bibliographystyle{apsr}
\bibliography{/home/erhard/Schreibtisch/Dropbox/Datenbanken/BibTex_files/lit_hist,/home/erhard/Schreibtisch/Dropbox/Datenbanken/BibTex_files/lit_mathsci,/home/erhard/Schreibtisch/Dropbox/Datenbanken/BibTex_files/lit_Weyl,/home/erhard/Schreibtisch/Dropbox/Datenbanken/BibTex_files/lit_scholz}

\end{document}